\documentclass[11pt]{amsart}

\usepackage{latexsym,amssymb,amsmath,graphics}
\usepackage[dvips]{graphicx}

\def\bb{{\mathcal B}}
\def\cc{{\mathcal C}}
\def\dd{{\mathcal D}}
\def\ee{{\mathcal E}}

\def\lll{{\mathcal L}}

\def\rr{{\mathcal R}}
\def\ss{{\mathcal S}}

\def\ffi{\varphi}

\def\dst{\displaystyle}

\def\C{{\mathbb{C}}}

\def\N{{\mathbb{N}}}
\def\P{{\mathbb{P}}}
\def\Q{{\mathbb{Q}}}
\def\R{{\mathbb{R}}}
\def\S{{\mathbb{S}}}

\def\X{{\mathbb{X}}}
\def\Y{{\mathbb{Y}}}

\def\fn{{\mathfrak{n}}}

\def\fN{{\mathfrak{N}}}
\def\fS{{\mathfrak{S}}}

\newcommand{\norm}[1]{{\left\|{#1}\right\|}}

\newcommand{\abs}[1]{{\left|{#1}\right|}}
\newcommand{\scal}[1]{{\left\langle{#1}\right\rangle}}

\newenvironment{notation}[1][]{\vskip1pt\noindent\rm\textit{Notation}\,:\ }{\rm\vskip1pt}
\newenvironment{remark}[1][]{\vskip1pt\noindent\rm\textit{Remark #1}\,:\ }{\rm\vskip1pt}
\newenvironment{definition}[1][]{\vskip3pt\noindent\sl\textbf{Definition.}\ }{\rm\vskip3pt}

\newtheorem{lemma}{Lemma}[section]
\newtheorem{proposition}[lemma]{Proposition}
\newtheorem{theorem}[lemma]{Theorem}
\newtheorem{corollary}[lemma]{Corollary}

\begin{document}

\title[Convolvable distributions on homogeneous Lie groups]{Distributions that are convolvable with generalized
Poisson kernel of solvable extensions of homogeneous Lie groups}
\author{Ewa Damek, Jacek Dziubanski, Philippe Jaming \& Salvador P\'erez-Esteva}
\address{E.D., J.D.\,: Institute of Mathematics\\ University of Wroclaw\\
50-384 Wroclaw\\
pl. Grunwaldzki 2/4\\ Poland}
\email{edamek@math.uni.wroc.pl, dziuban@math.uni.wroc.pl}

\address{Ph.J.\,: MAPMO-F\'ed\'eration Denis Poisson\\
Universit\'e d'Orl\'eans\\
BP 6759\\ F 45067 Orl\'eans Cedex 2\\
France}
\email{Philippe.Jaming@univ-orleans.fr}

\address{S.P.E.\,: Instituto de Matem\'{a}ticas, Unidad Cuernavaca\\Universidad Nacional Aut\'{o}noma de M\'{e}xico\\
Cuernavaca\\ Morelos 62251\\ M\'{e}xico}
\email{salvador@matcuer.unam.mx}

\begin{abstract}
In this paper, we characterize the class of distributions on an
homogeneous Lie group $\fN$ that can be extended via Poisson
integration to a solvable one-dimensional extension $\fS$ of $\fN$.
To do so, we introducte the $\ss'$-convolution on $\fN$
and show that the set of distributions that are $\ss'$-convolvable with
Poisson kernels is precisely the set of suitably weighted derivatives of
$L^1$-functions. Moreover, we show that the $\ss'$-convolution
of such a distribution with the Poisson kernel is harmonic and
has the expected boundary behaviour. Finally, we show that such
distributions satisfy some global weak-$L^1$ estimates.

\end{abstract}

\subjclass{}
\keywords{homogeneous Lie groups, distribution, $\ss'$-convolution,
Poisson integrals}

\date{September 16, 2004}

\thanks{Research partially financed by:\\ 
E.D., J.D., Ph.J.~:{\it European Commission}
Harmonic Analysis and Related Problems 2002-2006 IHP Network
(Contract Number: HPRN-CT-2001-00273 - HARP).\\
E.D., J.D.~:{\it European Commission Marie Curie Host Fellowship for the Transfer
of Knowledge} ``Harmonic Analysis, Nonlinear Analysis and
Probability'' MTKD-CT-2004-013389.\\
J.D~: Polish founds for science 2005--2008 (research project 1P03A03029).\\
S.P.-E.~: Conacyt-DAIC U48633-F.}

\maketitle

\tableofcontents

\section{Introduction}

The aim of this paper is to contribute to the understanding of the boundary behaviour
of harmonic functions on one dimensional extensions of homogeneous Lie groups. More precisely, we here address the question
of which distributions on the homogeneous Lie group can be extended via Poisson-like integration to the whole domain
and in which sense this distribution may be recovered as a limit on the boundary of its extension.
This question has been recently settled in the case of Euclidean harmonic functions on $\R^{n+1}_+$
in \cite{ddjpe.alvarez,ddjpe.AGPPE}. For sake of simplicity, let us detail the kind of results we are looking for in this context.

Let us endow $\R^{n+1}_+:=\{(x,t)\,:\ x\in\R^n, t>0\}$ with the Euclidean laplacian. The associated Poisson kernel
is then given by $\P_t(x)=\dst\frac{t}{(t^2+x^2)^{(n+1)/2}}$ and a compactly supported distribution $T$
can be extended into an harmonic function via convolution $u(x,t)=\P_t\ast T$. As $\P_t$ is not in the Schwartz class,
this operation is not valid for arbitrary distributions in $\ss'$. The question thus arizes of which
distributions in $\ss'$ can be extended via convolution with the Poisson kernel. The first task is to properly define
\emph{convolution} and it turns out that the best results are obtained by using the $\ss'$-convolution
which agrees with the usual convolution of distributions when this makes sense. The space of distributions
that can be $\ss'$-convolved with the Poisson kernel is then the space of derivatives of properly-weighted $L^1$-functions.
Moreover, the distribution obtained this way is a harmonic function which has the expected boundary behaviour.

In this paper, we generalize these results to one dimensional extensions of homogeneous Lie groups,
that is homogeneous Lie groups with a one-dimensional family of dilations acting on it. This is
a natural habitat for generalizing results on $\R^{n+1}_+$ and these spaces occur in various situations.
The most important to our sense is that homogeneous Lie groups occur in the Iwasawa decomposition of semi-simple
Lie groups and hence as boundaries of the associated rank one symmetric space
or more generally, as boundaries of homogeneous spaces of negative curvature
\cite{ddjpe.He}. Both symmetric spaces and homogeneous spaces of negavite
curvature are semi-direct products $\fS=\fN\bold R^*_+$ of a homogeneous group 
$\fN$ and $\bold R_+^*$ acting by dilations in the first case, or ``dilation like''
automorphisms in the second. For a large class of left-invariant operators
on $\fS$ bounded harmonic functions can be reproduced from their boundary
values  on $\fN$ via so called Poisson integrals. They involve Poisson
kernels whose behavior at infinity is very similar to the one of $\P_t$.
While for rank one symmetric spaces and the Laplace-Beltrami operator this
is immediate form an explicite formula, for the most general case it has
been obtained only recently after many years of considerable interest in 
the subject (see \cite{ddjpe.BDH} and references there).
Therefore, we consider a large family of kernels on which we only impose growth conditions that are similar
to those of usual Poisson kernels. This allows us to obtain the desired generalizations.

In doing so, the main
difficulty comes from the right choice of definition of the $\ss'$-convolution, since the various choices
are {\it a priori} non equivalent do to the non-commutative nature of the homogeneous Lie group.
Once the right choice is made, we obtain the full characterization of the space of distributions
the can be extended via Poisson integration. We then show that this extension has the desired properties, namely
that it is harmonic if the Poisson kernel is harmonic and that the original distribution is obtained as
a boundary value of its extension. Finally, we show that the harmonic functions obtained in this way
satisfy some global estimates.

The article is organized as follows. In the next section, we recall the main results on Lie groups that we will
use. We then devote a section to results on distributions on homogeneous Lie groups and the $\ss'$-convolution
on these groups. Section \ref{ddjpe.sec:main} is the main section of this paper. There we prove the characterization
of the space of distributions that are $\ss'$-convolvable with Poisson kernels and show that their $\ss'$-convolution
with the Poisson kernel has the expected properties. We conclude the paper by proving that functions
that are $\ss'$-convolution of distributions with the Poisson kernels satisfy global estimates.

\section{Background and preliminary results}

In this section we recall the main notations and results we need on homogeneous Lie algebras and groups.
Up to minor changes of notation, all results from this section that are given without proof can be found in
the first chapter of \cite{ddjpe.FS}, although in a different order.

\subsection{Homogeneous Lie algebras, norms and Lie groups}
\label{ddjpe.sec:homo}

Let $\fn$ be a real and finite dimensional nilpotent Lie algebra with Lie bracket denoted $[\cdot,\cdot]$.
We assume that $\fn$ is endowed with a family of dilations $\{\delta_a\,:\ a>0\}$, consisting of automorphisms of
$\fn$ of the form $\delta_a=\exp(A\log a)$ where $A$ is a diagonalizable linear operator on $\fn$ with positive 
eigenvalues. As usual, we will often write $a\eta$ for $\delta_a\eta$ and even $\eta/a$ for $\delta_{1/a}\eta$.
Without loss of generality, we assume that the smallest eigenvalue of $A$ is $1$.
We denote
$$
1=d_1\leq d_2\leq\cdots\leq d_n:=\bar d
$$
the eigenvalues of $A$ listed with multiplicity. We will write
$$
\Delta=\left\{\sum_{\alpha\in F}\alpha_jd_j
\,:\ F\subset\N^n\mbox{ finite}\right\}.
$$
If $\alpha$ is a multi-index, we will write $|\alpha|=\alpha_1+\cdots+\alpha_n$ for its length and
$d(\alpha)=d_1\alpha_1+\cdots+d_n\alpha_n$ for its weight.

Next, we fix a basis $X_1,\ldots,X_n$ of $\fn$ such that $AX_j=d_jX_j$ for each $j$
and write $\vartheta_1,\ldots,\vartheta_n$ for the dual basis of $\fn^*$. Finally we define an Euclidean structure on 
$\fn$ by declaring the $X_i$'s to be orthonormal.
The associated scalar product will be denoted $\scal{\cdot,\cdot}$ and the norm $\norm{\cdot}$.

We denote by $\fN$ the connected and simply connected Lie group that corresponds to $\fn$. If we denote by $V$
the underlying vector space of $\fn$ and by $\theta_k=\vartheta_k\circ\exp^{-1}$, then $\theta_1,\ldots,\theta_n$ form 
a system of global coordinates on $\fN$ that allow to see $\fN$ as $V$. Note that $\theta_k$ is homogeneous of degree 
$d_k$ in the sense that $\theta_k(\delta_a \eta)=a^{d_k}\theta_k(\eta)$.
The group law is then given by
$$
\theta_k(\eta\xi)=\theta_k(\eta)+\theta_k(\xi)+
\sum_{\alpha\not=0,\beta\not=0,d(\alpha)+d(\beta)=d_k}c_k^{\alpha,\beta}\theta^\alpha(\eta)\theta^\beta(\xi)
$$
for some constants $c_k^{\alpha,\beta}$ and $\theta^\alpha=\theta_i^{\alpha_1}\cdots\theta_n^{\alpha_n}$.
Note that the sum above only involves terms with degree of homogeneity $<d_k$, that is coordinates
$\theta_1,\ldots,\theta_{k-1}$. Although the group law is written multiplicatively, we will  write $0$ for the
identity of $\fN$.

Now we consider the semidirect products $\fS=\fN\rtimes\R ^*_+$
of such a nilpotent group $\fN$ with $\R ^*_+$,
that is, we consider $\fS=\fN\times\R^*_+$ with the multiplication
$$
(\eta , a)(\xi , b)=(\eta\delta _a(\xi),ab).
$$

Finally, we fix an homogeneous norm on $\fN$, that is a continous function $x\mapsto|x|$ from $\fN$ to $[0,+\infty)$
which is $\cc^\infty$ on $\fN\setminus\{0\}$ such that
\begin{enumerate}
\renewcommand{\theenumi}{\roman{enumi}}
\item $\abs{\delta_a\eta}=a\abs{\eta}$,

\item $\abs{\eta}=0$ if and only if $\eta=0$,

\item $\abs{\eta^{-1}}=\abs{\eta}$,

\item $\abs{\eta\cdot\xi}\leq(\abs{\eta}+\abs{\xi})$, $\gamma\geq 1$ and, according to
\cite{ddjpe.HS}, we may chose $|.|$ in such a way that $\gamma=1$,

\item this norm satisfies Petree's inequality: for $r\in\R$, 
$$
(1+|\eta\xi|)^r\leq (1+|\eta|)^{|r|}(1+|\xi|)^r.
$$
This inequality is obtained as follows: when $r\geq0$, write
$$
1+\abs{\xi\eta}\leq 1+(\abs{\eta}+\abs{\xi})\leq(1+\abs{\eta})(1+\abs{\xi})
$$
and raise it to the power $r$. For $r<0$, write
$$
1+\abs{\xi}\leq 1+(\abs{\xi\eta}+\abs{\eta^{-1}})\leq(1+\abs{\xi\eta}+\abs{\eta})
\leq(1+\abs{\xi\eta})(1+\abs{\eta})
$$
and raise it to the power $-r$.
\end{enumerate}
In particular, $d(\eta,\xi)=\abs{\eta^{-1}\xi}$ is a left-invariant metric on $\fN$.

For smoothness issues in the next sections, we will need the following notation.
Let $\Phi$ be a fixed $\C^\infty$ function on $[0,+\infty]$ such that $\Phi=1$ in $[0,1]$,
$\Phi(x)=x$ on $[2,+\infty)$ and $\Phi\geq 1$ on $[1,2]$. Then for $\mu\in\R$, we will denote by
$\omega_\mu(\eta)=(1+\Phi(|\eta|))^\mu$ which is $\cc^\infty$ in $\fN$.
In all estimates written bellow, $\omega_\mu$ can always be replaced by $(1+|\eta|)^\mu$.

\subsection{Haar measure and convolution of functions}

If $\eta\in\fN$ and $r>0$, we define
$$
B(\eta,r)=\{\xi\in\fN\,:\ |\xi^{-1}\eta|<r\}
$$
the ball of center $\eta$ and radius $r$. Note that $\overline{B(\eta,r)}$ is compact.

If $\mbox{d}\lambda$ denotes Lebesgue measure on $\fn$, then $\lambda\circ\exp^{-1}$ is a bi-invariant
Haar measure on $\fN$. We choose to normalize it so as to have $|B(\eta,1)|=1$ and still denote it by
$\mbox{d}\lambda$. Moreover, we have
$$
\abs{B(\eta,r)}=\abs{B(0,r)}=\abs{r\cdot B(0,1)}=r^Q,
$$
where $Q=d_1+\cdots+d_n=\mbox{tr}\,A$ is the homogeneous dimension of $\fN$.
This measure admits a polar decomposition. More precisely, if we denote by $\S=\{\eta\in\fN\,:\ |\eta|=1\}$,
there exists a measure $\mbox{d}\sigma$ on $\S$ such that for all $\ffi\in L^1(\fN)$,
$$
\int_{\fN}\ffi(\eta)\,\mbox{d}\lambda(\eta)=
\int_0^{+\infty}\int_{\S}\ffi(r\xi)r^{Q-1}\,\mbox{d}\sigma(\xi)\,\mbox{d}r.
$$
On $\fS$ the right-invariant Haar measure is given by $\dst\frac{\mbox{d}\lambda\,\mbox{d}a}{a}$.

Recall that the convolution on a group $\fN$ with left-invariant Haar measure $\mbox{d}\lambda$ is given by
$$
f\ast g(\eta)=\int_{\fN} f(\xi)g(\xi^{-1}\eta)\,\mbox{d}\lambda(\xi)=\int_G f(\eta\xi^{-1})g(\xi)\,\mbox{d}\lambda(\xi).
$$
This operation is not commutative but, writing $\check f(\eta)=f(\eta^{-1})$,
we have $f\ast g=(\check g\ast\check f)\check\ $.

We will need the following:

\begin{lemma}
Let $h$ be a $\cc^\infty$ function on $\fN$ supported in a
compact neighborhood of $0$ such that 
$$
\int_{\fN} h(\eta)\,\mathrm{d}\lambda(\eta)=1.
$$
Set $h_a(\eta)=a^{-Q} h(\delta_{a^{-1}} \eta)$, then the family $h_a$
forms a smooth compactly supported approximate identity.	In particular, if $f$ is continuous and bounded on $\fN$,
then $f\ast h_a\to f$ uniformly on compact sets as $a\to 0$.
\label{ddjpe.lem:approxid}
\end{lemma}

We will need the following elementary lemma that can be proved along the lines of \cite[Lemma 9]{ddjpe.AGPPE}:

\begin{lemma}
\label{ddjpe.lem:lem9}
For $r,s\in\R$, let 
$$
I_{r,s}(\eta)=\int_{\fN}(1+|\xi|)^r(1+|\xi^{-1}\eta|)^s
\,\mathrm{d}\lambda(\xi).
$$
Then, if $r+s+Q<0$, $I_{r,s}(\eta)$ is finite. Moreover, if this is the case, there is a constant $C_{r,s}$
such that, for every $\eta\in\fN$,
$$
I_{r,s}(\eta)\leq\begin{cases} C_{r,s}(1+\abs{\eta})^{r+s+Q}&\mbox{if }r+Q>0\mbox{ and }s+Q>0\\
C_{r,s}(1+\abs{\eta})^{\max(r,s)}\log(2+\abs{\eta})&\mbox{if }r+Q=0\mbox{ or }s+Q=0\\
C_{r,s}(1+\abs{\eta})^{\max(r,s)}&\mbox{else}
\end{cases}.
$$
\end{lemma}

\begin{proof} 
From Peetre's inequality
we immediately get the first part of the lemma.

From now on, we can assume that $r+s+Q<0$. Write $\fN=\Omega_1\cup\Omega_2\cup\Omega_3$
for a partition of $\fN$ given by
$$
\Omega_1=\left\{\xi\in\fN:\ \abs{\xi}\leq\frac{1}{2}\abs{\eta}\right\}
\quad\mbox{and}\quad
\Omega_2=\left\{\xi\in\fN:\ \abs{\xi}>\frac{1}{2}\abs{\eta},\ 
\abs{\xi^{-1}\eta}\leq\frac{1}{2}\abs{\eta}\right\}
$$
and let
$$
I_i(\eta)=\int_{\Omega_i}(1+\abs{\xi})^r(1+\abs{\xi^{-1}\eta})^s
\,\mbox{d}\lambda(\xi).
$$

First, for $\xi\in\Omega_1$, we have
$\frac{1}{2}\abs{\eta}\leq\abs{\xi^{-1}\eta}
\leq\frac{3}{2}\abs{\eta}$ so that
\begin{align}
I_1(\eta)\leq&
C_s(1+\abs{\eta})^s\int_{\Omega_1}(1+\abs{\xi})^r\,\mbox{d}\lambda(\xi)
\leq C_s(1+\abs{\eta})^s
\int_0^{\frac{\abs{\eta}}{2}}t^{Q-1}(1+t)^r\,\mbox{d}t\notag\\
\leq&\begin{cases}
C_{r,s}(1+\abs{\eta})^{r+s+Q}&\mbox{if }r+Q>0\\
C_{r,s}(1+\abs{\eta})^s\ln(2+\abs{\eta})&\mbox{if }r+Q=0\\
C_{r,s}(1+\abs{\eta})^s&\mbox{if }r+Q<0\\
\end{cases}.\notag
\end{align}

Next, for $\xi\in\Omega_2$, we have $\frac{1}{2}\abs{\eta}\leq\abs{\xi}\leq\frac{3}{2}\abs{\eta}$, thus
\begin{align}
I_2(\eta)\leq&
C_r(1+\abs{\eta})^r\int_{\Omega_2}(1+\abs{\xi^{-1}\eta})^s
\,\mbox{d}\lambda(\xi)
\leq C_r(1+\abs{\eta})^r
\int_0^{\abs{\eta}/2}t^{Q-1}(1+t)^s\,\mbox{d}t\notag\\
\leq&\begin{cases}
C_{r,s}(1+\abs{\eta})^{r+s+Q}&\mbox{if }s+Q>0\\
C_{r,s}(1+\abs{\eta})^r\ln(2+\abs{\eta})&\mbox{if }s+Q=0\\
C_{r,s}(1+\abs{\eta})^r&\mbox{if }s+Q<0\\
\end{cases}.\notag
\end{align}

Finally, for $\xi\in\Omega_3$, we have $\frac{1}{3}\abs{\xi}\leq\abs{\xi^{-1}\eta}\leq 3\abs{\xi}$
so that
\begin{align}
I_3(\eta)\leq&
C_{r,s}\int_{\Omega_3}(1+\abs{\xi})^r(1+\abs{\xi})^s
\mbox{d}\lambda(\xi)
\leq C_{r,s}\int_{\fN\setminus\Omega_1}
(1+\abs{\xi})^{r+s}\,\mbox{d}\lambda(\xi)\notag\\
=&C_{r,s}\int_{\frac{\abs{\eta}}{2}}^{+\infty} 
t^{Q-1}(1+t)^{r+s}\,\mbox{d}t
\leq C_{r,s}(1+\abs{\eta})^{r+s+Q}.\notag
\end{align}
The proof is then complete when grouping all estimates.
\end{proof}

\subsection{Invariant differential operators on $\fN$}
Recall that an element $X\in\fn$ can be identified with a left-invariant differential operator on $\fN$ via
$$
Xf(\xi)=\left.\frac{\partial}{\partial s}f\bigl(\xi.\exp(sX)\bigr)\right|_{s=0}.
$$
There is also a right-invariant differential operator $Y$ corresponding to $X$, given by
$$
Yf(\xi)=\left.\frac{\partial}{\partial s}f\bigl(\exp(sX).\xi\bigr)\right|_{s=0}.
$$
Note that $X$ and $Y$ agree at $\xi=0$. For $X_1,\ldots,X_n$ the basis of $\fn$ defined in section \ref{ddjpe.sec:homo}
we write $Y_1,\ldots,Y_n$ for the corresponding right-invariant differential operators.

If $\alpha$ is a multi-index, we will write
$$
\begin{matrix}
X^\alpha=X_1^{\alpha_1}\cdots X_n^{\alpha_n},&\widetilde{X}^\alpha=X_n^{\alpha_n}\cdots X_1^{\alpha_1},\\
Y^\alpha=Y_1^{\alpha_1}\cdots Y_n^{\alpha_n},&\widetilde{Y}^\alpha=Y_n^{\alpha_n}\cdots Y_1^{\alpha_1}.\\
\end{matrix}
$$
We will write $Z^\alpha$ if something is true for any of the above. For instance, we will use without further notice that
$$
|Z^\alpha\omega_\mu|\leq C\omega_{\mu-d(\alpha)}.
$$

For ``nice'' functions, one has\footnote{in \cite{ddjpe.FS} the
$\widetilde{\ }$ is missing, this is usually
harmless but not in this article.}
$$
\int_{\fN}X^\alpha f(\eta)g(\eta)\,\mbox{d}\lambda(\eta)
=(-1)^{|\alpha|}\int_{\fN}f(\eta)\widetilde{X}^\alpha g(\eta)
\,\mbox{d}\lambda(\eta)
$$
and
$$
\int_{\fN}Y^\alpha f(\eta)g(\eta)\,\mbox{d}\lambda(\eta)
=(-1)^{|\alpha|}\int_{\fN}f(\eta)\widetilde{Y}^\alpha g(\eta)
\,\mbox{d}\lambda(\eta).
$$
As a consequence, one also has
$$
\begin{matrix}
X^\alpha(f\ast g)=f\ast(X^\alpha g),&&\widetilde{X}^\alpha(f\ast g)=f\ast(\widetilde{X}^\alpha g),\\
Y^\alpha(f\ast g)=(Y^\alpha f)\ast g&\mbox{ and }&\widetilde{Y}^\alpha(f\ast g)=(\widetilde{Y}^\alpha g)\ast f.\\
\end{matrix}
$$
Moreover, using  $X^\alpha\check f=(-1)^{|\alpha|}(Y^\alpha f)\check{\ }$ or
$\widetilde{X}^\alpha\check f=(-1)^{|\alpha|}(\widetilde{Y}^\alpha f)\check{\ }$ and correcting
the proof in \cite{ddjpe.FS}, one gets
$$
(X^\alpha f)\ast g=f\ast(\widetilde{Y}^\alpha g) \mbox{ and } (\widetilde{X}^\alpha f)\ast g=f\ast(Y^\alpha g).
$$

Recall that a polynomial on $\fN$ is a function of the form
$$
P=\sum_{finite}a_\alpha\theta^\alpha
$$
and that its isotropic and homogeneous degrees are respectively defined by
$\max\{|\alpha|,\ a_\alpha\not=0\}$ and $\max\{d(\alpha),\ a_\alpha\not=0\}$.

For sake of simplicity,
we will write the Leibniz' Formula as
$$
X^\alpha(\ffi\psi)=\sum_{\beta\leq\alpha}\Lambda_{\alpha,\beta}X^\beta\ffi X^{\alpha-\beta}\psi,\ 
\widetilde{X}^\alpha(\ffi\psi)=\sum_{\beta\leq\alpha}\widetilde{\Lambda}_{\alpha,\beta}
\widetilde{X}^\beta\ffi \widetilde{X}^{\alpha-\beta}\psi.
$$

Further, we may write
\begin{equation}
\label{ddjpe.eq:jac1}
\widetilde{Y}^\alpha=\sum_{\beta\in\mathcal{I}_\alpha}
\widetilde{\mathbb{Q}}_{\alpha,\beta}X^\beta
\end{equation}
where $\mathcal{I}_\alpha=\{\beta\,:\ |\beta|\leq|\alpha|,
d(\beta)\geq d(\alpha)\}$ and
$\widetilde{\mathbb{Q}}_{\alpha,\beta}$ are homogeneous polynomials of homogeneous degree $d(\beta)-d(\alpha)$.
In the same way, any euclidean derivative can be written in terms of left or right invariant derivatives.
We will only need the following in the next section: for every $M$, there exist polynomials
$\omega_\alpha$, $|\alpha|\leq 2M$ and left-invariant operators $X^\alpha$ such that
\begin{equation}
\label{ddjpe.eq:jac1euc}
(I-\Delta)^M=\sum_{|\alpha|\leq 2M}\omega_\alpha X^\alpha.
\end{equation}

Finally, we will exhibit another link among several of this objects.
Let $h_a$ be as in Lemma \ref{ddjpe.lem:approxid} and let $f,\ffi$ be smooth compactly supported functions.
Then $\scal{(X^\alpha f)*h_a, \varphi}$ is
\begin{eqnarray*}
&=&\scal{X^\alpha f,\ffi*\check h_a}
=(-1)^{|\alpha|}\scal{f,\widetilde{X}^\alpha(\ffi*\check h_a)}
=\scal{f, \varphi*(\widetilde{Y}^\alpha h_a)^\vee}\\
&=&\int_{\fN}\int_{\fN} f(\xi)\ffi(\xi\eta)
(\widetilde{Y}^\alpha h_a)(\eta)
z,\mbox{d}\lambda(\eta)\,\mbox{d}\lambda(\xi)\\
&=&\int_{\fN}\int_{\fN}f(\xi)\ffi(\xi\eta)
\sum_{\beta\in\mathcal{I}_\alpha} 
\widetilde{\mathbb{Q}}_{\alpha,\beta}(\eta)(X^\beta h_a)(\eta)
\,\mbox{d}\lambda(\eta)\,\mbox{d}\lambda(\xi) \\
&=&\int_{\fN}\int_{\fN} f(\xi)
\sum_{\beta\in\mathcal{I}_\alpha} (-1)^{|\beta|}
\Big(\widetilde{X}^\beta 
\Big(\widetilde{\mathbb{Q}}_{\alpha,\beta}(\eta)
\ffi(\xi\eta)\Big)\Big) h_a(\eta)
\,\mbox{d}\lambda(\eta)\,\mbox{d}\lambda(\xi)\nonumber\\
&=&\int_{\fN}\int_{\fN} f(\xi)
\sum_{\beta\in\mathcal{I}_\alpha} (-1)^{|\beta|}
\sum_{\iota\leq\beta}\widetilde{\Lambda}_{\beta,\iota} 
\big(\widetilde{X}^{\beta-\iota}
\widetilde{\mathbb{Q}}_{\alpha,\beta}\big) (\eta)
\big(\widetilde{X}^\iota\ffi\big)(\xi\eta)h_a(\eta)
\,\mbox{d}\lambda(\eta)\,\mbox{d}\lambda(\xi).
\end{eqnarray*}
As $\widetilde{X}^{\beta-\iota}\widetilde{\mathbb{Q}}_{\alpha,\beta}$ is an homogeneous polynomial,
if it is not a constant, then $\widetilde{X}^{\beta-\iota}\widetilde{\mathbb{Q}}_{\alpha,\beta}(0)=0$. 
With Lemma \ref{ddjpe.lem:approxid}, it follows that
$$
\int_{\fN}\big(\widetilde{X}^{\beta-\iota}
\widetilde{\mathbb{Q}}_{\alpha,\beta}\big)(\eta)
\big(\widetilde{X}^\iota\ffi\big)(\xi\eta)h_a(\eta)
\,\mbox{d}\lambda(\eta)\to 0
$$
uniformly with respect to $\xi$ in compact sets, as $a\to 0$.
On the other hand, if $\widetilde{X}^{\beta-\iota}\widetilde{\mathbb{Q}}_{\alpha,\beta}$ is a constant,
\begin{eqnarray*}
\int_{\fN}\big(\widetilde{X}^{\beta-\iota}
\widetilde{\mathbb{Q}}_{\alpha,\beta}\big)(\eta)
\big(\widetilde{X}^\iota\ffi\big)(\xi\eta)h_a(\eta)\,\mbox{d}\lambda(\eta)
&=&\big(\widetilde{X}^{\beta-\iota}
\widetilde{\mathbb{Q}}_{\alpha,\beta}\big)(0)
\int_{\fN}\big(\widetilde{X}^\iota\ffi\big)(\xi\eta)h_a(\eta)
\,\mbox{d}\lambda(\eta)\\
&\dst\longrightarrow& \widetilde{X}^{\beta-\iota}\widetilde{\mathbb{Q}}_{\alpha,\beta}(0)\widetilde{X}^\iota\ffi(\xi)
\end{eqnarray*}
as $a\to0$, uniformly with respect to $\xi$ in compact sets, again with Lemma \ref{ddjpe.lem:approxid}.
We thus get that $\scal{(X^\alpha f)*h_a, \varphi}$ converges to
$$
\int_{\fN}f(\xi)\sum_{\beta\in\mathcal{I}_\alpha} (-1)^{|\beta|}
\sum_{\iota\leq\beta}\widetilde{\Lambda}_{\beta,\iota}
\big(\widetilde{X}^{\beta-\iota}\widetilde{\mathbb{Q}}_{\alpha,\beta}\big)(0)\widetilde{X}^\iota\ffi(\xi)\,\mbox{d}\lambda(\xi)
$$
On the other hand $(X^\alpha f)\ast h_a$ converges uniformly to $X^\alpha f$ on compact sets, thus
$$
\scal{(X^\alpha f)*h_a, \varphi}\to\scal{X^\alpha f,\ffi}=(-1)^{|\alpha|}\scal{f,\widetilde{X}^\alpha\ffi}.
$$
As the two forms of the limit are the same for all $f,\ffi$ with compact support, we thus get that
\begin{equation}
\label{ddjpe.eq:fundlink}
\widetilde{X}^\alpha= (-1)^{|\alpha|}\sum_{\beta\in\mathcal{I}_\alpha} (-1)^{|\beta|}
\sum_{\iota\leq\beta}\widetilde{\Lambda}_{\beta,\iota}
\big(\widetilde{X}^{\beta-\iota}\widetilde{\mathbb{Q}}_{\alpha,\beta}\big)(0)\widetilde{X}^\iota.
\end{equation}

\subsection{A decomposition of the Dirac distribution}

In Section \ref{ddjpe.sec:dpl1}, we will need the following result
about the existence of a parametrix:

\begin{lemma}
\label{ddjpe.lem:param}
For every integer $m$ and every compact set $K\subset\fN$
with $0$ in the interior, there exists 
and integer $M$, a family of left-invariant differential operators 
$X^\alpha$ of order $|\alpha|\leq M$ and a family of functions 
$\{F_\alpha\}_{|\alpha|\leq M}$ of class $\mathcal{C}^m$ with support
in $K$ such that
\begin{equation}
\label{ddjpe.eq:param}
\sum_{\alpha}X^\alpha F_\alpha=\delta_0
\end{equation}
where $\delta_0$ is the Dirac mass at origin.
\end{lemma}

\begin{proof}
Let us start with the proof in the euclidean case, eventhough this is
classical ({\it see} \cite{ddjpe.schwartz})

First, for $M$ big enough, the function $F_0$ defined on $\R^d$ by
$\widehat{F_0}(\xi)=1/(1+4\pi^2|\xi|^2)^M$ (where $\widehat{F}$ is the 
Fourier transform of $F$) is of class $\cc^m$ and satisfies
$(I-\Delta)^MF=\delta_0$ where $\Delta$ is the Euclidean Laplace 
operator.

Now let $\ffi$ be a smooth function supported in $K$ with $\ffi=1$
in a neighbourhood of $0$.
Then by Leibniz's rule, we get that $(I-\Delta)^M(F_0\ffi)$
is of the form
$$
\ffi(I-\Delta)^MF_0+\sum 
c_{\alpha\beta}\partial^{\beta}F_0\partial^\alpha\ffi.
$$
Note that $\partial^\alpha\ffi=0$ in a neighbourhood of $0$
and that $F_0$ is analytic away from $0$ so that, if we set
$\dst H=\sum_{0<|\alpha|\leq 2M,|\beta|\leq 2M}
c_{\alpha\beta}\partial^{\beta}F_0\partial^\alpha\ffi$
then $H$ is smooth and supported in $K$.
Further, as $(I-\Delta)^MF_0\ffi=\ffi(0)\delta_0=\delta_0$, we have
thus proved that there exists two functions $G$ of class $\mathcal{C}^m$ 
with support in $K$ such that
$$
(I-\Delta)^MG=\delta_0+H.
$$

To obtain \eqref{ddjpe.eq:param}, let us recall \eqref{ddjpe.eq:jac1euc}:
$$
(I-\Delta)^M=\sum_{|\alpha|\leq 2M}\omega_\alpha X^\alpha.
$$
It follows that, for $\psi\in\dd$,
\begin{eqnarray*}
\scal{(I-\Delta)^MG,\psi}
&=&\sum_{|\alpha|\leq 2M}
\scal{\omega_\alpha X^\alpha G,\psi}
=\sum_{|\alpha|\leq 2M}
(-1)^{|\alpha|}\scal{G,X^\alpha(\omega_\alpha\psi)}\\
&=&\sum_{|\alpha|\leq 2M}(-1)^{|\alpha|}
\sum_{\beta\leq\alpha}
\scal{G,X^{\alpha-\beta}\omega_\alpha X^\beta\psi}\\
&=&\sum_{|\alpha|\leq 2M}\sum_{\beta\leq\alpha}
\scal{X^\beta\bigl((-1)^{|\alpha|+|\beta|}G
X^{\alpha-\beta}\omega_\alpha\bigr),\psi}.
\end{eqnarray*}
We have thus written
$$
(I-\Delta)^M G=\sum_\alpha\sum_{\beta\leq\alpha}
X^\beta\bigl((-1)^{|\alpha|+|\beta|}GX^{\alpha-\beta}\omega_\alpha\bigr)
$$
which is of the desired form.
\end{proof}

\subsection{Laplace operators and Poisson kernels}
\label{ddjpe.sec:estpois}\ 

\begin{definition}
Let $\P$ be a smooth function on $\fN$ and let $\P_a(\eta)=a^Q\P(\delta_{a^{-1}}\eta)$.
We will say that $\P$ has property $(\mathcal{R}_\Gamma)$ if it satisfies the following estimates:
\begin{enumerate}
\renewcommand{\theenumi}{\roman{enumi}}
\item there exists a constant $C$ such that
$$
\frac{1}{C}\omega_{-Q-\Gamma}\leq \P \leq C\omega_{-Q-\Gamma};
$$
\item for every left-invariant operator $X^\alpha$, there is a constant $C_\alpha$ such that
for every $\eta\in\fN$, 
$|X^\alpha\P(\eta)|\leq C_\alpha \omega_{-Q-\Gamma-d(\alpha)}(\eta)$,

\item for every $k$, there is a constant $C_k$ such that
for every $\eta\in\fN$, and every $a>0$,
$$
|(a\partial_a)^k\P_a(\eta)|\leq C_\alpha a^{-Q}\omega_{-Q-\Gamma}(\delta_{a^{-1}}\eta).
$$
\end{enumerate}
\end{definition}

Note that several other important estimates will automatically result from this estimates.

\begin{enumerate}
\item First, by homogeneity of the left-invariant operator $X^\alpha$, there is a constant $C_\alpha$ such that
for every $\eta\in\fN$, and every $a>0$, 
$$
|X^\alpha\P_a(\eta)|\leq C_\alpha a^{-Q-d(\alpha)}\omega_{-Q-\Gamma-d(\alpha)}(\delta_{a^{-1}}\eta).
$$

\item Let $\X=X_{i_1}^{\alpha_1}\cdots X_{i_k}^{\alpha_k}$ be a left-invariant differential operator.
Set $d(\X)=d_{i_1}\alpha_1+\cdots+d_{i_k}\alpha_k$ its weight, then the commutation rules in $\fn$ imply that
$\X=\sum_{\beta:\ d(\beta)=d(\alpha)}c_\beta X^\beta$. It follows that
$$
|\X\P_a(\eta)|\leq C a^{-Q-d(\X)}\omega_{-Q-\Gamma-d(\X)}(\delta_{a^{-1}}\eta).
$$

\item Writing $Y^\alpha=\sum_{\beta\in\mathcal{I}_\alpha}\Q_{\alpha,\beta}X^\beta$ where $Q_{\alpha,\beta}$ is
an homogeneous polynomial of degree $d(\beta)-d(\alpha)$,
we get that
$$
|Y^\alpha\P_a(\eta)|\leq C a^{-Q-d(\alpha)}\omega_{-Q-\Gamma-d(\alpha)}(\delta_{a^{-1}}\eta).
$$
In particular, in all estimates, $\P_a$ can be replaced by $\check\P_a$. Also, as for the previous point,
$Y^\alpha$ may be replaced by $\Y=Y_{i_1}^{\alpha_1}\cdots Y_{i_k}^{\alpha_k}$.

\item The previous remark also shows that in point (ii) we may
as well impose the condition for \emph{right invariant} differential operators. This would not change the class of kernels.
\end{enumerate}

\begin{remark}
A large class of kernels satisfying property $\mathcal R_{\Gamma}$
is associated to left-invariant operators on $\fS$.
Consider a second order left-invariant operator on $\fS$ of the form
$$
\mathcal{L}=\sum _{j=1}^m Y_j^2 +Y.
$$
We assume the H\"ormander condition {\it i.e.} that
\begin{equation}\label{ddjpe.hor}
Y_1,...,Y_m \ \mbox{generate the Lie algebra of} \ \fS.
\end{equation}
The image of such an operator on $\R ^+$ under the natural
homomorphism $(\xi , a)\to a$ is, up to a multiplicative constant,
$$
(a\partial_a)^2-\alpha a\partial _a.
$$
If $a>0$ then there is a smooth integrable function $\P_a$ on $\fN$ such that the Poisson integrals
\begin{equation}
f*\P_a(\eta)=\int _\fN f(\xi)\P_a(\xi ^{-1}\eta )\ \mbox{d}\lambda(\xi)
\end{equation}
of $L^{\infty}$ function $f$ are $\mathcal{L}$-harmonic and moreover,
all bounded $\mathcal{L}$ harmonic functions are of this form. In particular, $\P _a(\eta)$ is $\mathcal{L}$-harmonic.

The properties i) and ii) for $\P $ have been proved in \cite{ddjpe.BDH} -
see the main theorem there for diagonal action and $\mathcal L$
satisfying \ref{ddjpe.hor}. iii) follows immediately from i) and the
(left-invariant) Harnack inequality applied to the harmonic function
$\P _a(\eta )$ i.e.
$$
|(a\partial _a)^k\P _a(\eta )|\leq C_k\P _a(\eta ).
$$
\end{remark}

Our first aim will be to give a meaning to such Poisson integrals
for as general as possible distributions $f$ so as to still obtain
an $\mathcal{L}$-harmonic functions when the kernel is $\mathcal{L}$-harmonic.

\section{Distributions on $\fN$}

\subsection{Basic facts and the space $\dd^{\prime}_{L^1}$}
\label{ddjpe.sec:dpl1}

Distributions on $\fN$ are defined as on $\R^n$ as the dual of the space $\dd:=\dd(\fN)$ of $\cc^\infty$
functions with compact support, endowed with the semi-norms given, for $K\subset\fN$ compact
and $\alpha$ a multi-index, by
$$
p_{K,\alpha}(\ffi)=\sup_{\eta\in K}\abs{\partial^\alpha(\ffi)}.
$$
We will write the space of distributions $\dd^{\prime}:=\dd^{\prime}(\fN)$,
endowed with the natural family of semi-norms. Notions such as support, Schwartz class $\ss:=\ss(\fN)$,
tempered distributions $\ss^\prime:=\ss^\prime(\fN)$,... are defined as for
distributions on $\R^n$ and the space of compactly supported distributions will be denoted
$\ee^{\prime}:=\ee^{\prime}(\fN)$.
Because of the link between left invariant derivatives and Euclidean derivatives
(similar to the links between left and right invariant derivatives, {\it see} \cite{ddjpe.FS}),
these spaces are just the usual spaces of derivatives on $\fN$ seen as $V\simeq \R^n$.
In particular, we will use the fact that every set of distributions that is weakly bounded is also
strongly bounded.

For $T\in\dd'$, we define $\check T\in\dd'$ by $\scal{\check T,\ffi}=\scal{T,\check\ffi}$,
while $X^\alpha T$ is defined by $\scal{X^\alpha T,\ffi}=(-1)^{|\alpha|}\scal{T,\widetilde{X}\ffi}$.

The definition of the convolution of two functions is easily extended to convolution
of a distribution with a smooth function via the following pairings:
for $T\in\dd'$ a distribution and $\psi,\ffi\in\dd$ smooth functions

--- the right convolution is given by $\scal{T\ast\psi,\ffi}=\scal{T,\ffi\ast\check\psi}$

--- the left convolution is given by $\scal{\psi\ast T,\ffi}=\scal{T,\check\psi\ast\ffi}$.

As in the Euclidean case, one may check that $T\ast\psi$ and $\psi\ast T$ are both smooth.

We will now introduce the space of integrable distributions $\dd^{\prime}_{L^1}$
and show that this is the space of derivatives of $L^1$ functions.

\begin{definition}
Let $\bb:=\bb(\fN)$ be the space of smooth functions $\ffi:\fN\to\C$ such that, for every 
left-invariant differential operator $X^\alpha$, $X^\alpha\ffi$ is bounded.

Let $\dot\bb:=\dot\bb(\fN)$ be the subspace of all $\ffi\in\bb(\fN)$ such that, for every 
left-invariant differential operator $X^\alpha$, $\abs{X^\alpha\ffi(u)}\to0$ when $\abs{u}\to\infty$.

We equip these spaces with the topology of uniform convergence of all derivatives.

The space $\dd^{\prime}_{L^1}=\dd^{\prime}_{L^1}(\fN)$ is the topological dual of $\dot\bb(\fN)$
endowed with the strong dual topology.
\end{definition}

Note that $\ss$ and $\cc^\infty_0$ are dense in $\dot\bb$ (but not in $\bb$) so that
$\dd^{\prime}_{L^1}$ is a subspace of $\ss'$. Note also that every compactly supported distribution is in
$\dd^{\prime}_{L^1}$.

It is also obvious that if $T\in\dd^{\prime}_{L^1}$, $\ffi\in\bb$ and $X^\alpha$ is left-invariant,
then $X^\alpha T\in\dd^{\prime}_{L^1}$ and $\ffi T\in\dd^{\prime}_{L^1}$.
We will need the following characterization of this space:
\begin{theorem} 
\label{ddjpe.th:d'l1}
Let $T\in\dd'(\fN)$. The following are equivalent
\begin{enumerate}
\renewcommand{\theenumi}{\roman{enumi}}
\item $T\in\dd^{\prime}_{L^1}(\fN)$;

\item $T$ has a representation of the form $T=\dst\sum_{finite}X^\alpha f_\alpha$
where $f_\alpha\in L^1(\fN)$ and $X^\alpha$ are left-invariant differential operators;

\item for every $\ffi\in\dd(\fN)$, the regularization $T\ast\ffi\in L^1(\fN)$.
\end{enumerate}
\end{theorem}

\begin{proof} The proof follows the main steps of the Euclidean case, see \cite[page 131]{ddjpe.schwartz}.
Denote by $\dd_1$ the set of all functions $\psi\in\dd$ such that $\norm{\psi}_\infty\leq1$.

\noindent$i)\Rightarrow iii)$ Assume that $T\in\dd^{\prime}_{L^1}$ and let $\ffi\in\dd$. 
 Now, note that
\begin{equation}
\label{ddjpe.eq:inter1}
\scal{T\ast\ffi,\psi}=\scal{T,\psi\ast\check\ffi}
\end{equation}
so, if $\ffi$ is fixed and $\psi$ runs over $\dd_1$, the set of numbers on the right of (\ref{ddjpe.eq:inter1})
is bounded, thus so is the set of numbers $\{\scal{T\ast\ffi,\psi},\ \psi\in\dd_1\}$. But $T\ast\ffi$ is a (smooth)
function so this implies that $T\ast\ffi\in L^1$.

\medskip

\noindent$iii)\Rightarrow ii)$ Assume that, for every $\psi\in\dd$, $T\ast\psi\in L^1$,
thus $T\ast\check\psi\in L^1$. 
Now, for $\psi\in\dd$ fixed, the set of numbers
$$
\scal{\check T\ast\ffi,\check\psi}=\scal{\check T,\check\psi\ast\check\ffi}=
\scal{T,\ffi\ast\psi}=\scal{T\ast\check\psi,\ffi}
$$
stays bounded when $\ffi$ runs over $\dd_1$. It follows that the set of distributions 
$\{\check T\ast\ffi,\ \ffi\in\dd_1\}$
is bounded in $\dd'$ since it is a weakly bounded set.

This implies that there exists an integer $m$ and a compact neighborhood $K$ of $0$ such that,
for every function $\psi$ of class $\cc^m$ with support in $K$, $\check T\ast\ffi\ast\psi(0)$ stays bounded
when $\ffi$ varies over $\dd_1$.
Using
$$
\check T\ast\ffi\ast\psi(0)=\scal{\check T\ast\ffi,\check\psi}=\scal{\check T,\check\psi\ast\check\ffi}=
\scal{T\ast\psi,\ffi}
$$
we get that $T\ast\psi\in L^1$ for every $\psi\in\cc^m$ with support in $K$.

\medskip

Now, according to Lemma \ref{ddjpe.lem:param}, we may write
$$
\sum_{finite}X^\alpha F_\alpha=\delta_0
$$
where the $F_\alpha$'s are of class $\cc^m$ and are supported in $K$.
It follows that
$$
T=\sum_{finite}T\ast\partial^\alpha F_\alpha=\sum_{finite}\partial^\alpha(T\ast F_\alpha).
$$
The first part of the proof shows that the $T\ast F_\alpha$'s are in $L^1$
so that we obtain the desired representation formula.

\noindent$ii)\Rightarrow i)$ is obvious so that the proof is complete.
\end{proof}

\begin{definition}
Let $\bb_c:=\bb_c(\fN)$ be the space $\bb(\fN)$ endowed with the topology for which $\ffi_n\to0$ if,

\begin{enumerate}
\renewcommand{\theenumi}{\roman{enumi}}
\item for every left-invariant differential operator $X^\alpha$, $X^\alpha\ffi_n\to0$
uniformly over compact sets,

\item for every left-invariant differential operator $X^\alpha$, the $X^\alpha\ffi_n$'s
are uniformly bounded.
\end{enumerate}
\end{definition}

The representation formula of $T\in\dd^{\prime}_{L^1}$ given by the previous theorem shows
that $T$ can be extended to a continuous linear functional on $\bb_c$ so that $\dd^{\prime}_{L^1}$ is also
the dual of $\bb_c$. For example, if we write
$T=f_0+\sum_{|\alpha|\geq 1} X^\alpha f_\alpha$, then
$$
\scal{T,1}_{\dd^{\prime}_{L^1},\bb_c}=\scal{f_0,1}
=\int_\fN f_0(\xi)\,\mbox{d}\lambda(\xi).
$$

\subsection{The $\ss'$-convoultion}

Recall that if $G\in\ss'$ and $\ffi\in\ss$ then $\check G\ast\ffi\in\cc^\infty$ so that the following definition makes sense:

\begin{definition}
Let $F,G\in\ss^{\prime}(\fN)$, we will say that they are $\ss^\prime$-convolvable if, for every $\ffi\in\ss(\fN)$,
$(\ffi\ast\check G)F\in\dd^{\prime}_{L^1}$. If this is the case, we define
$$
\scal{F\ast G,\ffi}
=\scal{(\ffi\ast\check G)F,1}_{\dd^{\prime}_{L^1},\bb_c}.
$$
\end{definition}

If $F,G\in\ss(\fN)$, then $F$ and $G$ are $\ss^\prime$-convolvable and the above definition coincides
with the usual one.
Indeed, for every $\ffi\in\ss(\fN)$,
\begin{eqnarray*}
\scal{F*G,\ffi}&=&\int_{\fN}F*G(\eta)\ffi(\eta)\,\mbox{d}\lambda(\eta)\\
&=&\int_{\fN}\int_{\fN}F(\xi)G(\xi^{-1}\eta)\ffi(\eta)
\,\mbox{d}\lambda(\xi)\,\mbox{d}\lambda(\eta)\\
&=&\int_{\fN}\left(\int_{\fN}
\ffi(\eta)\check G(\eta^{-1}\xi)\mbox{d}\lambda(\eta)\right)F(\xi).1
\,\mbox{d}\lambda(\xi)\\
&=&\scal{(\ffi*\check G)F,1}_{\dd^{\prime}_{L^1},\bb_c}.
\end{eqnarray*}

\begin{remark} There are various ways to define the $\ss^\prime$-convolution that extend the definition for functions.
For $S,T\in\dd'(\fN)$, let us cite the following:

\begin{enumerate}
\item $S$ and $T$ are $\ss^\prime_1$-convolvable if, for every $\ffi\in\dd(\fN)$,
$S_x\otimes T_y\ffi(xy)\in\dd^\prime_{L^1}(\fN\otimes\fN)$.
The $\ss^\prime_1$-convolution of $S$ and $T$ is then defined by 
$$
\scal{S\ast_1 T,\ffi}=
\scal{S_x\otimes T_y\ffi(xy),1}
_{\dd^\prime_{L^1}(\fN\otimes\fN),\bb_c(\fN\otimes\fN)}.
$$
\item $S$ and $T$ are $\ss^\prime_2$-convolvable if, for every $\ffi\in\dd$,
$S(\check T\ast\ffi)\in\dd^\prime_{L^1}(\fN)$ 
$$
\scal{S\ast_2 T,\ffi}=\scal{S(\check T\ast\ffi),1}
_{\dd^\prime_{L^1}(\fN),\bb_c(\fN)}.
$$
\item $S$ and $T$ are $\ss^\prime_3$-convolvable if, for every
$\ffi,\psi\in\ss(\fN)$, $(\check S\ast\ffi)(T\ast\check \psi)\in L^1(\fN)$.
The $\ss^\prime_3$-convolution of $S$ and $T$ is then defined by
$$
\scal{S\ast_3 T,\ffi\ast\psi}=\int_\fN
(\check S\ast\ffi)(\eta)( T\ast\check\psi)(\eta)\,\mbox{d}\lambda(\eta).
$$
\end{enumerate}

It turns out that in the Euclidean case, all four definitions are equivalent and lead to the same convolution
\cite{ddjpe.shiraishi}.
There are various obstructions to prove this in our situation,
mostly stemming from the fact that left and right-invariant derivatives
differ.

Also, one may replace the $\dd^\prime_{L^1}$ space by the similar
one defined with the help of right-invariant derivatives.
We will here stick to the choice given in the definition above
as it seems to us that this is the definition that gives
the most satisfactory results.
\end{remark}

One difficulty that arises is that the derivative of a convolution is not easily linked
to the convolution of a derivative. Here is an illustration of what may be done and of the difficulties that arise.
We hope that this will convince the reader that several facts that seem obvious 
(and are for usual convolutions of functions)
need to be proved, {\it e.g.} that $T\ast\P_a$ is harmonic if $\P_a$ is.

\begin{lemma} 
\label{ddjpe.lem:derconv}
Let $S,T\in\dd'(\fN)$ and
let $Y$ be a right-invariant differential operator of first order.
If $S$ and $T$ are $\ss^\prime$-convolvable, if $YS$ and $T$ are 
$\ss^\prime$-convolvable and if, for all $\ffi\in\ss(\fN)$, 
$Y\bigl((\ffi\ast\check T)S\bigr)\in\dd^\prime_{L^1}(\fN)$, then
$$
Y(S\ast T)=(Y S)\ast T.
$$
\end{lemma}

\begin{proof} As $(Y f)g=Y(fg)-fYg$, we get that
\begin{eqnarray*}
\scal{Y(S\ast T),\ffi}&=&-\scal{S\ast T,Y\ffi}=
-\scal{\bigl((Y\ffi)\ast\check T\bigr)S,1}\\
&=&-\scal{Y(\ffi\ast\check T)S,1}
=-\scal{Y\bigl((\ffi\ast\check T)S\bigr),1}
+\scal{(\ffi\ast\check T)YS,1}\\
&=&0+\scal{(YS)\ast T,\ffi}
\end{eqnarray*}
the next to last equality being justified by the assumptions on $F$, $G$.
\end{proof}

Using this lemma inductively gives
$$
Y^\alpha (S\ast T)=(Y^\alpha S)\ast T
$$
provided all intermediate steps satisfy the assumption of the lemma.
This is the case if $S$ is compactly supported.

\subsection{Weigthed spaces of distributions}
We will need the following weighted space of integrable distributions, introduced in the
Euclidean setting in \cite{ddjpe.horvath1,ddjpe.horvath2,ddjpe.ortner}.

\begin{definition}
\label{ddjpe.def:b}
Given $\mu \in \mathbb{R}$ we consider 
$$
\omega_\mu\dd_{L^{1}}^{\prime}(\fN)
:=\omega_\mu\dd_{L^{1}}^{\prime}(\fN)=
\left\{ T\in \ss^{\prime}(\fN):\ \omega_{-\mu }T\in \dd_{L^{1}}^{\prime }(\fN) \right\}
$$
with the topology induced by the map 
$$
\begin{matrix}
\omega_\mu\dd_{L^{1}}^{\prime}(\fN)& \rightarrow &\dd_{L^{1}}^{\prime }(\fN)\\
T& \mapsto &\omega_{-\mu }T
\end{matrix}.
$$
\end{definition}

This space admits an other representation given in the following lemma:

\begin{lemma}
\label{ddjpe.lem:tool1} 
Given $\mu\in\R$, we have
\begin{equation}
\omega_{\mu}\dd_{L^{1}}^{\prime}(\fN)=
\left\{ T\in S^{\prime}(\fN):\ T=\sum\limits_{finite}X^\alpha g_\alpha,
\ where\ g_\alpha\in L^1(\fN,\omega_{-\mu}d\lambda) \right\}.
\label{ddjpe.eq:weighted}
\end{equation}
\end{lemma}

\begin{proof}
Let us temporarily indicate with $\mathcal{V}$ the right hand side of (\ref{ddjpe.eq:weighted}).
Given $T\in \mathcal{V}$, we can write 
$\dst T=\sum\limits_{finite}X^\alpha\left(\omega_\mu f_\alpha\right)$, where $f_\alpha\in L^1$.
But then,
\begin{eqnarray*}
T&=&\sum\limits_{finite}\sum_{0\leq \beta \leq \alpha}\Lambda_{\alpha,\beta}
X^{\alpha-\beta}\omega_\mu X^\beta f_\alpha\\
&=&\omega_{\mu}\sum_{finite}\sum_{0\leq \beta \leq \alpha}\Lambda_{\alpha,\beta}
\omega_{-\mu}X^{\alpha-\beta}(\omega_{\mu}) X^\beta f_\alpha.
\end{eqnarray*}

By definition, the distribution $X^\beta f_\alpha$ belongs to 
$\dd_{L^{1}}^\prime$. Moreover, and easy computation shows that the function 
$\omega_{-\mu}X^{\alpha-\beta}(\omega_\mu)$ belongs to the
space $\bb$. Since $\dd_{L^{1}}^\prime$ is closed under multiplication by
functions in $\bb$, we conclude that $T$ belongs to
$\omega_\mu\dd_{L^{1}}^\prime$. 

Conversely, given $T\in\omega_\mu\dd_{L^{1}}^\prime$ we can write, by definition, 
$\dst T=\omega_\mu\sum_{finite}X^\alpha f_\alpha$, where $f_\alpha\in L^1$ or, 
$\dst T=\omega_\mu\sum_{finite}X^\alpha(\omega_{-\mu} g_\alpha)$, 
where $g_\alpha\in L^1( \omega_{-\mu}\mbox{d}\lambda)$. Now, given $\varphi\in\ss$,
the pairing $\scal{T,\varphi}_{\ss^\prime,\ss}$ can be written as 
$$
\sum\limits_{finite}(-1)^{|\alpha|}
\scal{g_\alpha,\omega_{-\mu}\widetilde{X}^\alpha\bigl(\omega_\mu\varphi\bigr)}_{\ss^{\prime },\ss}
=\sum\limits_{finite}\sum_{0\leq \beta \leq \alpha }
(-1)^{|\alpha|}\widetilde{\Lambda}_{\alpha,\beta }
\scal{g_\alpha,\omega_{-\mu}\left(\widetilde{X}^{\alpha-\beta}\omega_\mu\right) 
\widetilde{X}^\beta\varphi}_{\ss^{\prime },\ss}.
$$
We observe that for each multi-indexes $\alpha$ and $\beta$, the function 
$$
b_{\alpha,\beta}=(-1)^{|\alpha|}\widetilde{\Lambda}_{\alpha,\beta}\omega_{-\mu}
(\widetilde{X}^{\alpha -\beta }\omega_\mu)
$$
belongs to $\bb$. Thus,
$$
\scal{T,\varphi}_{\ss^{\prime },\ss}=\sum_{\alpha,\beta}(-1)^{|\beta|}
\scal{X^\beta\left((-1)^{|\beta|}b_{\alpha,\beta}g_\alpha\right),\varphi} _{\ss^{\prime },\ss}
$$
or,
$$
T=\sum_{\alpha ,\beta }(-1)^{|\beta|}X^{\beta }\left((-1)^{|\beta|}b_{\alpha,\beta }g_\alpha\right).
$$
To conclude that the distribution $T$ belongs to $\mathcal{V}$ we only need
to observe that $L^{1}\left(\omega_{-\mu}\mbox{d}\lambda\right)$ is closed under
multiplication by functions in $\bb$.
This completes the proof of Lemma \ref{ddjpe.lem:tool1}.
\end{proof}

As an immediate corollary, we get that
\begin{corollary}
\label{ddjpe.cor:strd'l1}
The space $\omega_{\mu}\dd^\prime_{L^1}(\fN)$ is closed under the action of left-invariant differential operators
$X^\alpha$ and under multiplication by functions in $\bb$.
\end{corollary}

\section{Distributions that are $\ss^\prime$-convolvable with the Poisson kernel}
\label{ddjpe.sec:main}

\subsection{Extensions of distributions with the Poisson kernel}

We are now in position to prove the following:

\begin{theorem}
\label{ddjpe.th:sconv}
Let $T\in\ss'$ and $\P$ be kernel with property $(\mathcal{R}_\Gamma)$. Then the following are equivalent:
\begin{enumerate}
\renewcommand{\theenumi}{\roman{enumi}}
\item $T\in\omega_{Q+\Gamma}\dd_{L^{1}}^{\prime }$,
\item $T$ is $\ss'$-convolvable with $\P_a$ for some $a>0$,
\item $T$ is $\ss'$-convolvable with $\P_a$ for each $a>0$.
\end{enumerate}
\end{theorem}

\begin{proof} It is of course enough to prove equivalence between (i) and (ii), the equivalence with (iii) will then automatically follow.

Let us assume that $T\in\omega_{Q+\Gamma}\dd_{L^{1}}^{\prime}$. We want to show that, if $\ffi\in\ss$,
$(\ffi\ast\check\P_a)T\in\dd_{L^1}^{\prime}$. It is enough to show that $(\ffi\ast\check\P_a)\omega_{Q+\Gamma}\in\bb$.
But, for a left-invariant derivative $X^\alpha$ and $\beta\leq\alpha$,
$$
X^{\alpha-\beta}(\ffi\ast\check\P_a)(\eta)
=\ffi\ast(X^{\alpha-\beta}\check\P_a)(\eta)
=\int_{\fN}X^{\alpha-\beta}\check\P_a(\eta\xi^{-1})\ffi(\xi)
\,\mbox{d}\lambda(\xi)$$
so that
\begin{align}
|X^{\alpha-\beta}(\ffi\ast\check\P_a)(\eta)|
\leq&C\int_{\fN}
\frac{1}{(1+|\eta\xi^{-1}|)^{Q+\Gamma+d(\alpha)-d(\beta)}}\,\ffi(\xi)
\,\mbox{d}\lambda(\xi)\notag\\
\leq&C\omega_{-Q-\Gamma-d(\alpha)+d(\beta)}(\eta)
\int_{\fN}(1+|\xi|)^{Q+\Gamma+d(\alpha)-d(\beta)}\ffi(\xi)
\,\mbox{d}\lambda(\xi)\notag
\end{align}
by Petree's inequality. As $|X^\beta\omega_{Q+\Gamma}|\leq C_\beta\omega_{Q+\Gamma-d(\beta)}$,
it follows from Leibnitz' Rule that $(\ffi\ast\check\P_a)\omega_{Q+\Gamma}\in\bb$.
The first part of the proof is thus complete.

Conversely, let us assume that $T$ is $\ss'$-convolvable with $P_a$ and
fix $\ffi\in\ss$, a non-negative function supported in $B(0,2)$ and such that $\ffi=1$
on $B(0,1)$. Then
\begin{align}
\ffi\ast\check \P_a(\eta)=&\int_{\fN}\ffi(\xi)\P_a(\eta^{-1}\xi)
\,\mbox{d}\lambda(\xi)\notag\\
\geq&C(a)\int_{B(0,1)}\frac{1}{\bigl(1+|\eta^{-1}\xi|\bigr)^{Q+\Gamma}}
\,\mbox{d}\lambda(\xi).\notag
\end{align}
But, for $\xi\in B(0,1)$, $|\eta^{-1}\xi|\leq (|\eta|+|\xi|)\leq (1+|\eta|$).
It follows that,
$$
\ffi\ast\check\P_a(\eta)\geq\frac{C(a)}{(1+|\eta|)^{Q+\Gamma}}\geq C(a)\omega_{-Q-\Gamma}(\eta).
$$
As we have already shown that $\omega_{Q+\Gamma}\ffi\ast\check\P_a(\eta)\in\bb$, we get that
$\frac{1}{\omega_{Q+\Gamma}\ffi\ast\check\P_a(\eta)}\in\bb$.
Finally, writing
$$
T=\omega_{Q+\Gamma}\frac{1}{\omega_{Q+\Gamma}
\ffi\ast\check\P_a(\eta)}(\ffi\ast\check \P_a)(\eta)T
$$
gives the desired representation since, by hypothesis, $(\ffi\ast\check\P_a)T\in\dd_{L^1}^{\prime}$.
\end{proof}

\subsection{Regularity of the $\ss^\prime$-convolution of a distribution and the Poisson kernel}
We may now prove the following lemma, which allows us to represent $T\ast P_a$ as a function:

\begin{lemma}
\label{ddjpe.lem:interp}
Let $T\in\omega_{Q+\Gamma}\dd^\prime_{L^1}(\fN)$ and $\P$ be a kernel with property $(\mathcal{R}_\Gamma)$. Then,
the $\ss^\prime$-convolution of $T$ with the kernel $\P_a$ is the function given by
\begin{equation}
\label{ddjpe.eq:inter}
\eta\mapsto\scal{\omega_{-Q-\Gamma}(\cdot)T,
\omega_{Q+\Gamma}(\cdot)\check\P_a(\eta^{-1}\cdot)}_{\dd^{\prime}_{L^1},\bb_c}
\end{equation}
\end{lemma}

\begin{proof} First note that $\xi\mapsto \omega_{Q+\Gamma}(\xi)\check\P_a(\eta^{-1}\xi)$ is
in $\bb$ and $\omega_{-Q-\Gamma}(\cdot)T\in\dd^\prime_{L^1}$ so that (\ref{ddjpe.eq:inter}) makes sense.

We want to prove that, if $\dst T=\sum_{finite}\omega_{Q+\Gamma}X^\alpha f_\alpha$ with $f_\alpha\in L^1$,
and if $\ffi\in\ss$, then $\scal{T\ast\P_a,\ffi}:=\scal{(\ffi\ast\check\P_a)T,1}$ is
equal to 
$\scal{\scal{\omega_{-Q-\Gamma}(\cdot)T,
\omega_{Q+\Gamma}(\cdot)\check\P_a(\eta^{-1}\cdot)}_{\dd^{\prime}_{L^1},\bb_c},\ffi(\eta)}$.

By linearity, it is enough to consider only one term in the sum, $T=\omega_{Q+\Gamma}X^\alpha f$
with $f\in L^1(\fN)$. But then
\begin{eqnarray*}
\scal{\omega_{Q+\Gamma}(\ffi\ast\check\P_a)X^\alpha f,1} &=&
(-1)^{|\alpha|}\scal{f,\widetilde{X}^\alpha
\bigl(\omega_{Q+\Gamma}(\ffi\ast\check\P_a)\bigr)}\\
&=&(-1)^{|\alpha|}
\sum_{\beta\leq\alpha}\widetilde{\Lambda}_{\alpha,\beta}\int_{\fN}
f(\xi)\widetilde{X}^{\alpha-\beta}\omega_{Q+\Gamma}(\xi)
\widetilde{X}^\beta(\ffi\ast\check\P_a)(\xi)
\,\mbox{d}\lambda(\xi).
\end{eqnarray*}
Further, we have
$$
\widetilde{X}^\beta(\ffi\ast\check\P_a)
=\ffi\ast(\widetilde{X}^\beta\check\P_a)(\xi)
=\int_{\fN}\ffi(\eta)(\widetilde{X}^\beta\check\P_a)(\eta^{-1}\xi)
\,\mbox{d}\lambda(\eta).
$$
It follows that
$\scal{\omega_{Q+\Gamma}(\ffi\ast\check\P_a)X^\alpha f_\alpha,1}$ is 
\begin{eqnarray}
&=&
(-1)^{|\alpha|}\sum_{\beta\leq\alpha}\widetilde{\Lambda}_{\alpha,\beta}
\int_{\fN}\int_{\fN}
f(\xi)\widetilde{X}^{\alpha-\beta}
\omega_{Q+\Gamma}(\xi)(\widetilde{X}^\beta\check\P_a)(\eta^{-1}\xi)
\,\mbox{d}\lambda(\xi)\,\ffi(\eta)\,\mbox{d}\lambda(\eta)
\label{ddjpe.eq:invder}\\
&=&(-1)^{|\alpha|}\int_{\fN}\left(\int_{\fN}f(\xi)
\widetilde{X}^\alpha\bigl(\omega_{Q+\Gamma}(\cdot)
\check\P_a(\eta^{-1}\cdot)\bigr)(\xi)
\,\mbox{d}\lambda(\xi)\right)\ffi(\eta)\,\mbox{d}\lambda(\eta)
\nonumber
\end{eqnarray}
using $(\widetilde{X}^\beta\check\P_a)(\eta^{-1}\xi)=\widetilde{X}^\beta_\xi\check\P_a(\eta^{-1}\xi)$ and
Leibnitz' Rule. Thus
$$
\scal{\omega_{Q+\Gamma}(\ffi\ast\check\P_a)X^\alpha f,1}
=\int_{\fN}\scal{X^\alpha f(\xi),\omega_{Q+\Gamma}(\xi)\check\P_a(\eta^{-1}\xi)}
\,\ffi(\eta)\,\mbox{d}\lambda(\eta)
$$
as claimed.

All intervertions of integrals are easily justified by the fact that $\omega_{Q+\Gamma}\check\P_a\in\bb$.
\end{proof}

\begin{corollary}
\label{ddjpe.cor:reg}
Let $T\in\omega_{Q+\Gamma}\dd^\prime_{L^1}$ and $\P$ be a kernel with property $(\mathcal{R}_\Gamma)$.
Then the function $T\ast \P_a$ is smooth. Moreover,
for any left-invariant derivative $X^\alpha$, $T$ is $\ss^\prime$-convolvable with $X^\alpha\P_a$ and
$X^\alpha(T\ast\P_a)=T\ast(X^\alpha\P_a)$ 
and for any $k\in\N$, $T$ is $\ss^\prime$-convolvable with $(a\partial_a)^k\P_a$ and
$(a\partial_a)^k(T\ast\P_a)=T\ast\bigl((a\partial_a)^k\P_a\bigr)$.
In particular, $T\ast\P_a$ is harmonic if $\P$ is.
\end{corollary}

\begin{proof}
As the proof of the implication $i)\Rightarrow ii)$ of Theorem \ref{ddjpe.th:sconv} only depends
on the estimates of the Poisson kernel
from Section \ref{ddjpe.sec:estpois}, we get with the same proof that if $T\in\omega_{Q+\Gamma}\dd^{\prime}_{L^1}$
then $T$ is $\ss^\prime$-convolvable with $X^\alpha\P_a$ and $(a\partial_a)^k\P_a$.

For the other assertions, we may again assume that
$T=\omega_{Q+\Gamma}X^\alpha f$. As $T\ast\P_a$ is a function,
from \eqref{ddjpe.eq:invder} in the proof of the previous lemma, we get that
$$
T\ast\P_a(\eta)=(-1)^{|\alpha|}
\sum_{\beta\leq\alpha}\widetilde{\Lambda}_{\alpha,\beta}\int_{\fN}f(\xi)
\widetilde{X}^{\alpha-\beta}\omega_{Q+\Gamma}(\xi)
\widetilde{X}_\xi^\beta\bigl(\check\P_a(\eta^{-1}\xi)\bigr)
\,\mbox{d}\lambda(\xi).
$$
It then remains to differente with respect to $\eta$ under the integral.
\end{proof}

We will need the space $\dd_{L^1(\omega_\mu)}$ of all functions $\ffi\in\cc^\infty$
such that, for every left-invariant partial differential operator $X^\alpha$,
$X^\alpha\ffi\in L^1(\omega_\mu\mbox{d}\lambda)$ endowed with the topology given by the familly of 
semi-norms
$$
\norm{\ffi}_{\alpha,\mu}=\sum_{\beta\leq\alpha}\norm{X^\beta\ffi}_{L^1(\omega_\mu\mbox{d}\lambda)}.
$$

We may get a more precise estimate of the Poisson integrals at fixed level.
\begin{proposition}
\label{ddjpe.prop:reg3}
Let $T\in\omega_{Q+\Gamma}\dd^\prime_{L^1}(\fN)$. For each $a>0$, the $\ss^\prime$-convolution $T\ast \P_a$ belongs to
$\dd_{L^1(\omega_{-Q-\Gamma}\mathrm{d}\lambda)}$.
\end{proposition}

\begin{proof} By linearity, it is enough to prove that, if $T=\omega_{Q+\Gamma}X^\alpha f$
for some $f\in L^1(\fN)$, then $X^\iota(T\ast \P_a)=T\ast X^\iota\P_a\in \omega_{Q+\Gamma}\dd_{L^1}$.
But, from (\ref{ddjpe.eq:inter}), we get that
\begin{eqnarray}
T\ast X^\iota\P_a(\eta)
&=&(-1)^{|\alpha|}\scal{f,\widetilde{X}^\alpha\bigl(\omega_{Q+\Gamma}(\cdot)X^\iota\check\P_a(\eta^{-1}\cdot)\bigr)}
\nonumber\\
&=&(-1)^{|\alpha|}
\sum_{\beta\leq\alpha}\widetilde{\Lambda}_{\alpha,\beta}
\int_{\fN}f(\xi)\widetilde{X}^{\alpha-\beta}\omega_{Q+\Gamma}(\xi)
\widetilde{X}^\beta X^\iota\check\P_a(\eta^{-1}\xi)
\,\mbox{d}\lambda(\xi)
\label{ddjpe.eq:prop36}
\end{eqnarray}
using Leibnitz' Formula and the facts that $f\in L^1$ and $\omega_{Q+\Gamma}(\cdot)X^\iota\P_a(\eta^{-1}\cdot)\in\bb$.
Using the estimates
$$
|\widetilde{X}^{\alpha-\beta}\omega_{Q+\Gamma}(\xi)|\leq C(\alpha,\beta)\omega_{Q+\Gamma-d(\alpha)+d(\beta)}(\xi)
$$
and
$$
|\widetilde{X}^\beta_\xi X^\iota_\xi\check\P_a(\eta^{-1}\xi)|
\leq C(\alpha,\beta,a)\omega_{-Q-\Gamma-d(\beta)-d(\iota)}(\eta^{-1}\xi)
$$
we see that the $L^1(\omega_{-Q-\Gamma}\mbox{d}\lambda)$-norm of each term of the sum in (\ref{ddjpe.eq:prop36}) is bounded by
\begin{align}
C\int_{\fN}&\omega_{-Q-\Gamma}(\eta)\int_{\fN}|f(\xi)|
\omega_{Q+\Gamma-d(\alpha)+d(\beta)}(\xi)
\omega_{-Q-\Gamma-d(\beta)-d(\iota)}(\eta^{-1}\xi)
\,\mbox{d}\lambda(\xi)\,\mbox{d}\lambda(\eta)\notag\\
&=\,C\int_{\fN}|f(\xi)|\omega_{Q+\Gamma-d(\alpha)+d(\beta)}(\xi)
\int_{\fN}\omega_{-Q-\Gamma}(\eta)
\omega_{-Q-\Gamma-d(\beta)-d(\iota)}(\eta^{-1}\xi)
\,\mbox{d}\lambda(\eta)\,\mbox{d}\lambda(\xi)\notag\\
&\leq\,C\int_{\fN} |f(\xi)|
\omega_{-d(\alpha)+d(\beta)}(\xi)\,\mbox{d}\lambda(\xi)\notag
\end{align}
with Lemma \ref{ddjpe.lem:lem9}. As $d(\beta)\leq d(\alpha)$
we get the desired result.
\end{proof}

\subsection{The Dirichlet problem in $\omega_{Q+\Gamma}\dd^\prime_{L^1}$}

We will now prove that $T$ is the boundary value of $T\ast\P_a$
in the $\omega_{Q+\Gamma}\dd^\prime_{L^1}$
sense. 

\begin{theorem}
\label{ddjpe.th:dirichlet2}
Let $T\in\omega_{Q+\Gamma}\dd^\prime_{L^1}$, then the convolution $T\ast \P_a$ converges to $T$
in $\omega_{Q+\Gamma}\dd^\prime_{L^1}$ when $a\to 0^+$.
\end{theorem}

\begin{proof}
We want to prove that, for $\ffi\in\dot\bb$,
\begin{equation}
\scal{\omega_{-Q-\Gamma}(T\ast\P_a),\ffi}_{\dd^{\prime}_{L^1},\dot\bb}
\to
\scal{\omega_{-Q-\Gamma}T,\ffi}_{\dd^{\prime}_{L^1},\dot\bb}
\label{ddjpe.eq:converge2}
\end{equation}
when $a\to0$.
It is of course enough to consider $T=\omega_{Q+\Gamma}X^\alpha f$ with $f\in L^1$.

Write $\ffi_{-Q-\Gamma}=\omega_{-Q-\Gamma}\ffi$ and ${}_\xi\ffi_{-Q-\Gamma}(\eta)=\ffi_{-Q-\Gamma}(\xi\eta)$.
Then $\scal{\omega_{-Q-\Gamma}(T\ast\P_a),\ffi}_{\dd^{\prime}_{L^1},\dot\bb}$ is 
\begin{eqnarray*}
&=&(-1)^{|\alpha|}\int_{\fN}\int_{\fN}f(\xi)\widetilde{X}^\alpha_\xi
\bigl(\omega_{Q+\Gamma}(\xi)\check\P_a(\eta^{-1}\xi)\bigr)
\ffi_{-Q-\Gamma}(\eta)\,\mbox{d}\lambda(\eta)\,\mbox{d}\lambda(\xi)\\
&=&(-1)^{|\alpha|}\int_{\fN}\int_{\fN}f(\xi)
\sum_{\beta\leq\alpha}\widetilde{\Lambda}_{\alpha,\beta}
(\widetilde{X}^{\alpha-\beta}\omega_{Q+\Gamma})(\xi)
\widetilde{X}^\beta_\xi\bigl(\check\P_a(\eta^{-1}\xi)\bigr)
\ffi_{-Q-\Gamma}(\eta)\,\mbox{d}\lambda(\eta)\,\mbox{d}\lambda(\xi)\\
&=&(-1)^{|\alpha|}\int_{\fN}\int_{\fN}f(\xi)
\sum_{\beta\leq\alpha}\widetilde{\Lambda}_{\alpha,\beta}
(\widetilde{X}^{\alpha-\beta}\omega_{Q+\Gamma})(\xi)
(\widetilde{X}^\beta\check\P_a)(\eta^{-1})
\ffi_{-Q-\Gamma}(\xi\eta)\,\mbox{d}\lambda(\eta)\,\mbox{d}\lambda(\xi)\\
&=&(-1)^{|\alpha|}\int_{\fN}\int_{\fN}f(\xi)
\sum_{\beta\leq\alpha}\widetilde{\Lambda}_{\alpha,\beta}
(\widetilde{X}^{\alpha-\beta}\omega_{Q+\Gamma})(\xi)
(-1)^{|\beta|}\widetilde{Y}^\beta\P_a(\eta){}_\xi\ffi_{-Q-\Gamma}(\eta)
\,\mbox{d}\lambda(\eta)\,\mbox{d}\lambda(\xi).
\end{eqnarray*}
Now let $\psi$ be a smooth cut-off function such that $\psi(\eta)=1$
if $|\eta|\leq 1$ and $\psi(\eta)=0$
if $|\eta|\geq 2$ and write $\tilde\psi=1-\psi$. Then
$\scal{\omega_{-Q-\Gamma}(T\ast\P_a),\ffi}_{\dd^{\prime}_{L^1},\dot\bb}=S_1+S_2$
where $S_1$ is
$$
(-1)^{|\alpha|}\int_{\fN}\int_{\fN}f(\xi)
\sum_{\beta\leq\alpha}\widetilde{\Lambda}_{\alpha,\beta}
(\widetilde{X}^{\alpha-\beta}\omega_{Q+\Gamma})(\xi)
(-1)^{|\beta|}\widetilde{Y}^\beta\P_a(\eta)
{}_\xi\ffi_{-Q-\Gamma}(\eta)\psi(\eta)
\,\mbox{d}\lambda(\eta)\,\mbox{d}\lambda(\xi)
$$
while $\dst S_2=(-1)^{|\alpha|}\sum_{\beta\leq\alpha}
\widetilde{\Lambda}_{\alpha,\beta}S_2^\beta$ with
$$
S_2^\beta=\int_{\fN}\int_{\fN}f(\xi)
(\widetilde{X}^{\alpha-\beta}\omega_{Q+\Gamma})(\xi)
(-1)^{|\beta|}\widetilde{Y}^\beta\P_a(\eta)
{}_\xi\ffi_{-Q-\Gamma}(\eta)\tilde \psi(\eta)
\,\mbox{d}\lambda(\eta)\,\mbox{d}\lambda(\xi).
$$

\medskip

Let us first show that $S_2\to 0$, that is each $S_2^\beta\to0$.
As, for $|\eta|\geq 1$,
$$
(1+|\eta|/a)^{-Q-\Gamma-d(\beta)}\leq a^{Q+\Gamma+d(\beta)}|\eta|^{-Q-\Gamma-d(\beta)}\leq
Ca^{Q+\Gamma+d(\beta)}(1+|\eta|)^{-Q-\Gamma-d(\beta)},
$$
so that using the estimates of derivatives of $\P_a$ and $\omega_{Q+\Gamma}$,
we get that
\begin{eqnarray*}
|S_2^\beta|&\leq&C\int_{\fN}|f(\xi)|
(1+|\xi|)^{Q+\Gamma-d(\alpha)+d(\beta)}\times\\
&&\qquad\qquad\times\int_{|\eta|\geq 1}
a^{Q+d(\beta)}(1+|\eta|/a)^{-Q-\Gamma-d(\beta)}(1+|\xi\eta|)^{-Q-\Gamma}
\,\mbox{d}\lambda(\eta)\,\mbox{d}\lambda(\xi)\\
&\leq&Ca^\Gamma
\int_{\fN}|f(\xi)|(1+|\xi|)^{Q+\Gamma-d(\alpha)+d(\beta)}\times\\
&&\qquad\qquad\times\int_{|\eta|\geq 1}
(1+|\eta|)^{-Q-\Gamma-d(\beta)}(1+|\xi\eta|)^{-Q-\Gamma}
\,\mbox{d}\lambda(\eta)\,\mbox{d}\lambda(\xi)\\
&\leq&Ca^\Gamma\norm{f}_{L^1}
\end{eqnarray*}
with Lemma \ref{ddjpe.lem:lem9}. It follows that $S_2\to 0$.

\medskip

Let us now turn to $S_1$. First, from \eqref{ddjpe.eq:jac1},
\begin{eqnarray*}
S_1&=&(-1)^{|\alpha|}\int_{\fN}\int_{\fN}f(\xi)\sum_{\beta\leq\alpha}
\widetilde{\Lambda}_{\alpha,\beta}
(\widetilde{X}^{\alpha-\beta}\omega_{Q+\Gamma})(\xi)
\times\\
&&\qquad\ \times(-1)^{|\beta|}\sum_{\iota\in\mathcal{I}_\beta}
\widetilde{\mathbb{Q}}_{\beta,\iota}(\eta)X^\iota\P_a(\eta)
\,{}_\xi\ffi_{-Q-\Gamma}(\eta)\psi(\eta)
\,\mbox{d}\lambda(\eta)\,\mbox{d}\lambda(\xi)\\
&=&(-1)^{|\alpha|}\int_{\fN}\int_{\fN}f(\xi)
\sum_{\beta\leq\alpha}\widetilde{\Lambda}_{\alpha,\beta}
(\widetilde{X}^{\alpha-\beta}\omega_{Q+\Gamma})(\xi)
\times\\
&&\qquad\ \times\ 
(-1)^{|\beta|}\P_a(\eta)\sum_{\iota\in\mathcal{I}_\beta}(-1)^{|\iota|}
\widetilde{X}^\iota\bigl(\widetilde{\mathbb{Q}}_{\beta,\iota}
\,{}_\xi\ffi_{-Q-\Gamma}\psi\bigr)(\eta)
\,\mbox{d}\lambda(\eta)\,\mbox{d}\lambda(\xi)\\
&=&(-1)^{|\alpha|}\int_{\fN}\int_{\fN}f(\xi)
\sum_{\beta\leq\alpha}\widetilde{\Lambda}_{\alpha,\beta}
(\widetilde{X}^{\alpha-\beta}\omega_{Q+\Gamma})(\xi)
\times\\
&&\qquad\ \times\ (-1)^{|\beta|}\P_a(\eta)\sum_{\iota\in\mathcal{I}_\beta}(-1)^{|\iota|}
\sum_{\iota'\leq\iota}\widetilde{\Lambda}_{\iota,\iota'}
\widetilde{X}^{\iota-\iota'}\bigl(\widetilde{\mathbb{Q}}_{\beta,\iota}\,{}_\xi\ffi_{-Q-\Gamma}\bigr)(\eta)
\widetilde{X}^{\iota'}\psi(\eta)\,\mbox{d}\lambda(\eta)\,\mbox{d}\lambda(\xi).
\end{eqnarray*}

Now, if $\iota'\not=0$, $\widetilde{X}^{\iota'}\psi$ is supported in $1\leq|\eta|\leq 2$.
Further, from Leibnitz' Rule, $\ffi\in\bb$ and Peetre's inequality we get that
$\widetilde{X}^{\iota-\iota'}(\widetilde{\mathbb{Q}}_{\beta,\iota}\,{}_\xi\ffi_{-Q-\Gamma})(\eta)$ is bounded by $C\omega_{-Q-\Gamma}(\xi)$
with $C$ independent from $\eta$. It follows that
$$
\abs{\int_{\fN}\P_a(\eta)\widetilde{X}^{\iota-\iota'}
\bigl(\widetilde{\mathbb{Q}}_{\beta,\iota}\,{}_\xi
\ffi_{-Q-\Gamma}\bigr)(\eta)\widetilde{X}^{\iota'}\psi(\eta)
\,\mbox{d}\lambda(\eta)}
\leq C\omega_{-Q-\Gamma}(\xi)
\int_{1\leq|\eta|\leq 2}\P_a(\eta)\,\mbox{d}\lambda(\eta).
$$
Consequently, since this integral goes to $0$, we have
$$
\int_{\fN}\int_{\fN}f(\xi)(\widetilde{X}^{\alpha-\beta}
\omega_{Q+\Gamma})(\xi)\P_a(\eta)\widetilde{X}^{\iota-\iota'}
\bigl(\widetilde{\mathbb{Q}}_{\beta,\iota}\,{}_\xi
\ffi_{-Q-\Gamma}\bigr)(\eta)\widetilde{X}^{\iota'}\psi(\eta)
\,\mbox{d}\lambda(\eta)\,\mbox{d}\lambda(\xi)\to 0.
$$
This shows that $S_1$ has same limit as
\begin{eqnarray*}
S_1^\prime&:=&(-1)^{|\alpha|}\int_{\fN}\int_{\fN}f(\xi)
\sum_{\beta\leq\alpha}\widetilde{\Lambda}_{\alpha,\beta}
(\widetilde{X}^{\alpha-\beta}\omega_{Q+\Gamma})(\xi)
\times\\
&&\qquad\ \times\ 
(-1)^{|\beta|}\P_a(\eta)\sum_{\iota\in\mathcal{I}_\beta}(-1)^{|\iota|}
\widetilde{X}^\iota\bigl(\widetilde{\mathbb{Q}}_{\beta,\iota}\,{}_\xi
\ffi_{-Q-\Gamma}\bigr)(\eta)
\psi(\eta)\,\mbox{d}\lambda(\eta)\,\mbox{d}\lambda(\xi)\\
&=&(-1)^{|\alpha|}\int_{\fN}\int_{\fN}f(\xi)
\sum_{\beta\leq\alpha}\widetilde{\Lambda}_{\alpha,\beta}
(\widetilde{X}^{\alpha-\beta}\omega_{Q+\Gamma})(\xi)\times\\
&&\qquad\ \times\ (-1)^{|\beta|}\P_a(\eta)\sum_{\iota\in\mathcal{I}_\beta}(-1)^{|\iota|}
\sum_{\iota'\leq\iota}\widetilde{\Lambda}_{\iota,\iota'}
\widetilde{X}^{\iota-\iota'}\widetilde{\mathbb{Q}}_{\beta,\iota}(\eta)\widetilde{X}^\iota{}_\xi\ffi_{-Q-\Gamma}(\eta)
\psi(\eta)\,\mbox{d}\lambda(\eta)\,\mbox{d}\lambda(\xi).
\end{eqnarray*}
Now, if
$\widetilde{X}^{\iota-\iota'}\widetilde{\mathbb{Q}}_{\beta,\iota}$
is not a constant polynomial, then
$\widetilde{X}^{\iota-\iota'}\widetilde{\mathbb{Q}}_{\beta,\iota}(0)=0$ 
so that 
\begin{equation}
\label{ddjpe.eq:phil1}
\int_{\fN}\P_a(\eta)\widetilde{X}^{\iota-\iota'}
\widetilde{\mathbb{Q}}_{\beta,\iota}(\eta)
\widetilde{X}^\iota{}_\xi\ffi_{-Q-\Gamma}(\eta)\psi(\eta)
\,\mbox{d}\lambda(\eta)
\end{equation}
goes to $0$ when $a\to 0$, while if $\widetilde{X}^{\iota-\iota'}\widetilde{\mathbb{Q}}_{\beta,\iota}$ is constant, then, as $\psi(0)=1$, this integral goes to
$$
\widetilde{X}^{\iota-\iota'}\widetilde{\mathbb{Q}}_{\beta,\iota}(0)
\widetilde{X}^\iota{}_\xi\ffi_{-Q-\Gamma}(0).
$$
Moreover, as \eqref{ddjpe.eq:phil1} stays bounded by $C\omega_{-Q-\Gamma}(\xi)$,
from the Dominated convergence theorem,
we get that
\begin{eqnarray*}
S_1^{\prime}&\to&
(-1)^{|\alpha|}\int_{\fN}f(\xi)
\sum_{\beta\leq\alpha}\widetilde{\Lambda}_{\alpha,\beta}
(\widetilde{X}^{\alpha-\beta}\omega_{Q+\Gamma})(\xi)\times\\
&&\qquad\qquad\times\ 
(-1)^{|\beta|}\sum_{\iota\in\mathcal{I}_\beta}(-1)^{|\iota|}
\sum_{\iota'\leq\iota}\widetilde{\Lambda}_{\iota,\iota'}
\widetilde{X}^{\iota-\iota'}
\widetilde{\mathbb{Q}}_{\beta,\iota}(0)
\widetilde{X}^\iota{}_\xi\ffi_{-Q-\Gamma}(0)\,\mbox{d}\lambda(\xi)\\
&=&(-1)^{|\alpha|}\int_{\fN}f(\xi)
\sum_{\beta\leq\alpha}\widetilde{\Lambda}_{\alpha,\beta}
(\widetilde{X}^{\alpha-\beta}\omega_{Q+\Gamma})(\xi)
\widetilde{X}^\beta{}_\xi\ffi_{-Q-\Gamma}(0)
\,\mbox{d}\lambda(\xi)
\end{eqnarray*}
where we have used Identity (\ref{ddjpe.eq:fundlink}) in the last equality.
But $\widetilde{X}^\beta{}_\xi\ffi_{-Q-\Gamma}(0)=\widetilde{X}^\beta\ffi_{-Q-\Gamma}(\xi)$
so that Leibnitz' Formula implies that this limit is
$$
(-1)^{|\alpha|}\int_{\fN}f(\xi)\widetilde{X}^\alpha
\bigl(\omega_{Q+\Gamma}\ffi_{-Q-\Gamma}\bigr)(\xi)
\,\mbox{d}\lambda(\xi)
=(-1)^{|\alpha|}\int_{\fN}f(\xi)\widetilde{X}^\alpha\ffi(\xi)
\,\mbox{d}\lambda(\xi)=\scal{X^\alpha f,\ffi}
$$
as claimed.
\end{proof}

\begin{remark}
Assume that $\P_a$ is harmonic for some left-invariant differential operator $\lll$.
The above result imply that given for a distribution $T\in
\omega_{Q+\Gamma}\dd_{L^1}^\prime$, the function $u=T\ast\P_a$ is a solution of
the Dirichlet problem 
$$
\left\{\begin{array}{lc}
\lll u=0 & \text{in}\ \fS \\ 
u|_{a=0}=T &  
\end{array}\right.
$$
where the boundary condition is now interpreted in the sense of convergence in $\omega_{Q+\Gamma}\dd^\prime_{L^1}$
as $a\rightarrow 0^{+}$.
\end{remark}

\section{Global estimates for Poisson integrals of distributions in $\omega_{Q+\Gamma}\dd^\prime_{L^1}$}

In this section, we will prove that the Poisson integrals of distributions in $\omega_{Q+\Gamma}\dd^\prime_{L^1}$
satisfy some global integrability conditions. In order to do so, we will use the following notations.

\begin{notation}
For a Borel set $F\subset \fS$, we denote by $|F|$ its measure 
with respect to $\mbox{d}\lambda\,\mbox{d}a$. A function on $\fS$
is said to be in $L^{1,\infty}(\mbox{d}\lambda\,\mbox{d}a)$
if there exists a constant $C$ such that, for all $\alpha>0$,
$$
\abs{\{(\eta,a)\in\fS\,:\ |f(\eta,a)|>\alpha\}}\leq\frac{C}{\alpha}.
$$
For $\Gamma\geq1$, let $\dst\Phi_\Gamma(\eta,a)=\frac{a^\Gamma}{(a+|\eta|)^{Q+\Gamma}}$
and note that
$\dst\frac{1}{a}\Phi_\Gamma\in L^{1,\infty}(\mbox{d}\lambda\,\mbox{d}a)$.
Indeed
\begin{eqnarray*}
|\{(\eta,a)\,:\ \frac{1}{a}\Phi_\Gamma(\eta,a)>\alpha\}|
&=&\int_0^{\alpha^{-1/(Q+1)}}
|B(0,a^{(\Gamma-1)/(Q+\Gamma)}\alpha^{-1/(Q+\Gamma)}-a)|\,\mbox{d}a\\
&=&\int_0^{\alpha^{-1/(Q+1)}}
\left((a^{(\Gamma-1)/(Q+\Gamma)}\alpha^{-1/(Q+\Gamma)}-a\right)^Q
\,\mbox{d}a\\
&=&\frac{1}{\alpha}\int_0^1\left(t^{(\Gamma-1)/(Q+\Gamma)}-t\right)^Q
\,\mbox{d}t
\end{eqnarray*}
by changing variable $t=a\alpha^{1/(Q+1)}$. It should also be noted that
$\dst\Phi_\Gamma\notin 
L^{1,\infty}\left(\frac{\mbox{d}\lambda\,\mbox{d}a}{a}\right)$.

We will denote by $\mathcal{M}_\Gamma$ the set of measures $\mu$ on $\fN$ such that
$$
\int_{\fN}(1+|\xi|)^{-(Q+\Gamma)}\,\mbox{d}|\mu|(\xi)<+\infty.
$$
For $\mu\in\mathcal{M}_\Gamma$ and $\eta\in\fN$, let us denote by $\mu_\eta$ the left translate of $\mu$ by $\eta$,
that is the measure defined by
$$
\int_{\fN}\ffi(\xi)\,\mbox{d}\mu_\eta(\xi)
=\int_\fN\ffi(\eta\xi)\,\mbox{d}\mu(\xi)
$$
for all continuous functions $\ffi$ with compact support on $\fN$. 
From Petree's inequality, we get that $\mu_\eta\in \mathcal{M}_\Gamma$.
Further note that, if $\mu\in\mathcal{M}_\Gamma$, then
$$
\mu\bigl(B(0,r)\bigr)\leq (1+r)^{Q+\Gamma}\int_{|\xi|<r}(1+|\xi|)^{-(Q+\Gamma)}\,\mbox{d}\mu(\xi)
\leq C(1+r)^{Q+\Gamma}.
$$
\end{notation}

We are now in position to prove the following:

\begin{theorem} 
Let $\Gamma\geq1$ and let $\P$ be a kernel with property $(\rr_\Gamma)$.
If $\mu\in\mathcal{M}_\Gamma$ then
\begin{equation}
\label{ddjpe.eq:thsj1}
\frac{1}{a}(1+a+|\eta|)^{-Q-\Gamma}\mu\ast\P_a(\eta)
\in L^{1,\infty}(\mathrm{d}\lambda\,\mathrm{d}a).
\end{equation}
Moreover, for every $a_0>0$,
\begin{equation}
\label{ddjpe.eq:sjogren2}
(1+a+|\eta|)^{-Q-\Gamma}a^{-\Gamma}\mu\ast\P_a(\eta)
\chi_{\{(\eta,a)\in\fS\,: a>a_0\}}(\eta,a)\in 
L^{\infty}(\mathrm{d}\lambda\,\mathrm{d}a).
\end{equation}
\label{ddjpe.th:sjogren}
\end{theorem}
\begin{proof}
Let $E_l=\{(\eta,a)\in\fS\:\ |\eta|\leq 1\mbox{ and }a\leq 1\}$ and $E_g=\fS\setminus E_l$
and let $\eta_0\in\fN$ be such that $|\eta_0|\geq 2$.
Assume that we have proved that for every $E\subset E_g$
$$
\frac{1}{a}(1+a+|\eta|)^{-Q-\Gamma}\mu\ast\P_a(\eta)
\chi_E\in L^{1,\infty}(\mbox{d}\lambda\,\mbox{d}a)
$$
for every measure $\mu\in\mathcal{M}_\Gamma$.
Without loss of generality we may assume that $\mu$ is a positive measure
Applying this to $E=\eta_0^{-1}E_l$ and to the left-translate
$\mu_0\in\mathcal{M}_\Gamma$ of $\mu$ by $\eta_0^{-1}$ we get that
\begin{align}
(1+a+|\eta|)^{-Q-\Gamma}a^{-1}&\mu\ast\P_a(\eta)\chi_{E_l}(\eta,a)\notag\\
&=(1+a+|\eta|)^{-Q-\Gamma}a^{-1}\mu_0\ast\P_a(\eta_0^{-1}\eta)\chi_{\eta_0^{-1}E_l}(\eta_0^{-1}\eta,a)\notag\\
&\leq C(1+a+|\eta_0^{-1}\eta|)^{-Q-\Gamma}a^{-1}
\mu_0\ast\P_a(\eta_0^{-1}\eta)\chi_{\eta_0^{-1}E_l}(\eta_0^{-1}\eta,a)\notag\\
&\in L^{1,\infty}(\mbox{d}\lambda\mbox{d}a).\notag
\end{align}

It is thus enough to prove that
\begin{equation}
\label{ddjpe.eq:thsj2}
\frac{1}{a}(1+a+|\eta|)^{-Q-\Gamma}\mu\ast\P_a(\eta)
\chi_{E_g}\in L^{1,\infty}(\mbox{d}\lambda\mbox{d}a).
\end{equation}

Note that if $\P$ has property $(\rr_\Gamma)$ then
$(1+a+|\eta|)^{-Q-\Gamma}a^{-1}\mu\ast\P_a(\eta)$ is bounded by
\begin{eqnarray*}
&&C\frac{a^{(\Gamma-1)}}{(1+a+|\eta|)^{Q+\Gamma}}\int_{\fN}\frac{\mbox{d}\mu(\xi)}{(a+|\eta^{-1}\xi|)^{Q+\Gamma}}\\
&=&C\frac{a^{(\Gamma-1)}}{(1+a+|\eta|)^{Q+\Gamma}}\left(\int_{|\xi|\leq\frac{1}{2}|\eta|}
+\int_{\frac{1}{2}|\eta|\leq|\xi|\leq 2|\eta|}
+\int_{2|\eta|\leq|\xi|}\right)\frac{\mbox{d}\mu(\xi)}{(a+|\eta^{-1}\xi|)^{Q+\Gamma}}\\
&=&I+II+III.
\end{eqnarray*}

Let us first estimate $I$. Note that, if $|\xi|\leq\frac{1}{2}|\eta|$, then
$$
a+|\eta^{-1}\xi|\geq a+|\eta|-|\xi|\geq a+\frac{1}{2}|\eta|\geq C(a+|\eta|)
\geq (1+|\eta|)/2
$$
since we only consider $(\eta,a)\in E_g$. It follows that
\begin{eqnarray*}
I\chi_{E_g}&\leq&C\frac{a^{(\Gamma-1)}}{(1+a+|\eta|)^{Q+\Gamma}}
\int_{|\xi|\leq\frac{1}{2}|\eta|}\mbox{d}\mu(\xi)(1+|\eta|)^{-Q-\Gamma}\chi_{E_g}\\
&\leq&C\frac{1}{a}\frac{a^\Gamma}{(1+a+|\eta|)^{Q+\Gamma}}\in L^{1,\infty}(\mbox{d}\lambda\mbox{d}a).
\end{eqnarray*}
Moreover, this computation also shows that $a^{-\Gamma+1}I\in L^\infty(\mbox{d}\lambda\mbox{d}a)$.

\medskip

Let us now estimate $III$. Note that, if $|\xi|\geq 2|\eta|$, then
\begin{equation}
\label{ddjpe.eq:estIII}
a+|\eta^{-1}\xi|\geq a+|\xi|-|\eta|\geq a+\frac{1}{2}|\xi|\geq C(a+|\xi|).
\end{equation}
Further, as $(\eta,a)\in E_g$ then either $a\geq 1$ or $|\eta|\geq 1$ in which case $|\xi|\geq 2$. Therefore
$a+|\eta^{-1}\xi|\geq C(1+|\xi|)$. It follows that
$$
III\chi_{E_g}\leq C\frac{a^{(\Gamma-1)}}{(1+a+|\eta|)^{Q+\Gamma}}\int_{\fN}(1+|\xi|)^{-Q-\Gamma}\mbox{d}\mu(\xi)
\leq C\frac{1}{a}\frac{a^\Gamma}{(1+a+|\eta|)^{Q+\Gamma}}\in L^{1,\infty}(\mbox{d}\lambda\mbox{d}a).
$$
Again, the same computation shows that $a^{-\Gamma+1}III\in L^\infty(\mbox{d}\lambda\mbox{d}a)$.

\medskip

We will now prove the result for $II$. To do so, notice first that, if $a\geq a_0$ then
$$
II\leq\frac{C}{a^{Q+1}}\int_{|\xi|\leq2|\eta|}d\mu(\xi)\leq
C\frac{1}{a}\frac{a^\Gamma}{a^{Q+\Gamma}(1+|\eta|)^{Q+\Gamma}}
\leq C a^{-Q-1}
$$
and again $a^{-\Gamma+1}II\chi_{a>a_0}\in L^\infty(\mbox{d}\lambda\mbox{d}a)$.
It is thus enough to prove the result for $II\chi_{a\leq a_0}$.

Now note that, if $\frac{1}{2}|\eta|\leq|\xi|\leq 2|\eta|$,
then $(1+a+|\eta|)^{-Q-\Gamma}\leq C(1+|\xi|)^{-Q-\Gamma}$, thus
$$
II\leq C\frac{1}{a}\int_{\frac{1}{2}|\eta|\leq|\xi|\leq 2|\eta|}\Phi_\Gamma(\eta^{-1}\xi,a)
(1+|\xi|)^{-Q-\Gamma}\mbox{d}\mu(\xi)
=C\frac{1}{a}\int_{\frac{1}{2}|\eta|\leq|\xi|\leq 2|\eta|}\Phi_\Gamma(\eta^{-1}\xi,a)
\mbox{d}\nu(\xi)
$$
where $\nu$ is a finite measure on $\fN$. Thus $II$ is estimated with the help of the following proposition:

\begin{proposition}
For every finite positive measure $\nu$ on $\fN$, the function $U_\nu$ defined on $\fS$ by
$$
U_\nu(\eta,a)=\int_{\fN}
\frac{a^{\Gamma-1}}{(a+|\eta^{-1}\xi|)^{Q+\Gamma}}\,\mathrm{d}\nu(\xi)
$$
belongs to $\dst L^{1,\infty}(\mathrm{d}\lambda\,\mathrm{d}a)$.
\label{ddjpe.prop:sjogren}
\end{proposition}

The proof will follow a simplified version of that of Theorem 1 in \cite{ddjpe.Sj} which deals with
the Euclidean case, for more general measures.

\begin{notation}
On $\fS$ we denote by $D_\infty$ the distance given by
$D_\infty\bigl((\eta,a),(\eta',a')\bigr)=
\max\bigl(|\eta^{-1}\eta'|,|a-a'|\bigr)$.
\end{notation}

\begin{proof}[Proof of Proposition \ref{ddjpe.prop:sjogren}]
We want to prove that
$U_\nu\in L^{1,\infty}(\mbox{d}\lambda\,\mbox{d}a)$, that is,
that there exists a constant $C\geq0$ such that for all $\alpha>0$,
$$
|\{(\eta,a)\in\fS\,:\ U_\nu(\eta,a)>\alpha\}|\leq\frac{C}{\alpha}.
$$
For $i_0$ a non-negative integer, let
$$
K_0=B(0,2^{i_0})\times]0,2^{i_0}].
$$
It is enough to prove that there is a constant $C\geq0$, independent of $i_0$ such that, for all $\alpha>0$
$$
|\{(\eta,a)\in\fS\,:\ U_\nu(\eta,a)>\alpha\}\cap K_0|\leq\frac{C}{\alpha}.
$$
To do so, we will show that there is a constant $C$
such that, for each $\alpha>0$, we may construct a set $S$ 
which satisfies the following properties:
\begin{enumerate}
\renewcommand{\theenumi}{\roman{enumi}}
\item\label{ddjpe.pointi} $|\{(\eta,a)\,:\ U_\nu>\alpha\}\cap K_0|\leq C|S|$;

\item\label{ddjpe.pointii} $U_\nu(\eta,a)>\frac{\alpha}{C}$ for all $(\eta,a)\in S$;

\item\label{ddjpe.pointiii} for all $\eta\in\fN$,
$$
U_S(\eta):=\int_S\frac{a^{\Gamma-1}}{(a+|\eta^{-1}\xi|)^{Q+\Gamma}}
\,\mbox{d}\lambda(\xi)\,\mbox{d}a
$$
satisfies $U_S(\eta)\leq C$.
\end{enumerate}

Once this is done, we can conclude as follows
\begin{eqnarray*}
|\{(\eta,a)\in\fS\,:\ U_\nu(\eta,a)>\alpha\}\cap K_0|&\leq&
C|S|\ \leq\ \frac{C^2}{\alpha}\int_S U_\nu(\eta,a)
\,\mbox{d}\lambda(\eta)\,\mbox{d}a\\
&=&\frac{C^2}{\alpha}\int_{\fN}U_S(\eta)\,\mbox{d}\nu(\eta)\\
&\leq&\frac{C^3}{\alpha}\norm{\nu}
\end{eqnarray*}
where we have respectively used Property \ref{ddjpe.pointi}, \ref{ddjpe.pointii}, Fubini's Theorem and Property
\ref{ddjpe.pointiii}.

\medskip

\noindent{\sl Construction of the set $S$.}

\smallskip

We will use a dyadic covering of $K_0$:

\noindent--- Set $Q_i=B(0,2^{i_0})\times[2^{i_0-i-1},2^{i_0-i}]$, $i=0,1,\ldots$

\noindent--- cover each $Q_i$ by sets of the form
$$
Q_{i,j}=B(\eta_{i,j},2^{i_0-i-3})\times[2^{i_0-i-1},2^{i_0-i}]
$$
in such a way that each element of $Q_i$
belongs to at most $\kappa$ sets $Q_{i,j}$ where $\kappa$ is a set
that depends only on the group $\fN$. This is possible thanks to
a covering lemma that may be found {\it e.g.} in \cite[Section 1.F]{ddjpe.FS}.

We will order the $Q_{i,j}$'s by lexicographic order and define inductively
the authorized pieces $A_{i,j}$ and the associated set of forbidden pieces $F_{i,j}$ as follows:

\noindent--- $A_{i,j}=Q_{i,j}$ if
\begin{enumerate}
\renewcommand{\theenumi}{\alph{enumi}}
\item $|Q_{i,j}\cap\{(\eta,a)\in K_0\,:\ U_\nu(\eta,a)>\alpha\}|>0$,
\item and $Q_{i,j}\notin\dst\bigcup\limits_{(l,k)<(i,j)}F_{l,k}$.
\end{enumerate}
Else, we set $A_{i,j}=\emptyset$.

\noindent--- if $A_{i,j}\not=\emptyset$ we define the set of forbidden pieces as
$$
F_{i,j}=\{Q_{l,k}\,:\ (l,k)>(i,j)\mbox{ and }
D_\infty(Q_{l,k},Q_{i,j})<2^{i_0-i+\frac{l-i}{Q+\Gamma}+1}\}.
$$
Else we set $F_{i,j}=\emptyset$.

Note that the authorized pieces are disjoint and that, if $A_{i,j}\not=\emptyset$, then $F_{i,j}$ has the following property:
\begin{equation}
\label{ddjpe.eq:fijaij}
\abs{\bigcup\limits_{Q\in F_{i,j}}Q}\leq C|A_{i,j}|.
\end{equation}

\begin{proof}[Proof of (\ref{ddjpe.eq:fijaij})] Assume that
$Q_{l,k}\in F_{i,j}$ and let $(\eta,a)\in Q_{l,k}$. Then
\begin{eqnarray*}
\mathrm{d}(\eta,\eta_{i,j})&\leq& \mathrm{d}\bigl(\eta,\fN\setminus B(\eta,2^{i_0-l-3})\bigr)+
D_\infty\bigl(Q_{l,k},Q_{i,j}\bigr)
+\mathrm{d}\bigl(\eta_{i,j},\fN\setminus B(\eta_{i,j},2^{i_0-i-3})\bigr)\\
&<&2^{i_0-l-3}+2^{i_0-i+\frac{l-i}{Q+\Gamma}+1}+2^{i_0-i-3}\\
&\leq&2^{i_0-i+\frac{l-i}{Q+\Gamma}+2}.
\end{eqnarray*}
It follows that
$$
Q_{l,k}\subset B(\eta_{i,j},C2^{i_0-i+\frac{l-i}{Q+\Gamma}})\times[2^{i_0-l-1},2^{i_0-l}].
$$
Now, as pieces of different order are disjoint,
$$
\abs{\bigcup\limits_{Q\in F_{i,j}}Q}=\sum_{m=0}^{+\infty}
\abs{\bigcup\limits_{Q_{i+m,k}\in F_{i,j}}Q_{i+m,k}}
$$
and as those of a given order overlap at most $\kappa$ times, this is
\begin{eqnarray*}
&\leq&\kappa\sum_{m=0}^{+\infty}\abs{B(\eta_{i,j}
,C_22^{i_0-i+\frac{m}{Q+\Gamma}+1})\times[2^{i_0-i-m-1},2^{i_0-i-m}]}\\
&\leq&C\sum_{m=0}^{+\infty}2^{(i_0-i+\frac{m}{Q+\Gamma})Q}2^{i_0-i-m}\\
&=&C2^{(i_0-i)(Q+1)}\sum_{m=0}^{+\infty}2^{-\frac{\Gamma}{Q+\Gamma}m}
=C|A_{i,j}|.
\end{eqnarray*}
\end{proof}

Finally, we set $S=\dst\bigcup_{(i,j)}A_{i,j}$.

\medskip

\noindent{\sl Proof of Property \ref{ddjpe.pointi}.}

\smallskip

By construction, the authorized and the forbidden pieces cover 
$\{U_\nu>\lambda\}\cap K_0$ and as these overlap at most $\kappa$
times, we obtain
\begin{eqnarray*}
|\{(\eta,a)\in K_0\,: U_\nu(\eta,a)>\lambda\}|&\leq&
\sum_{(i,j)}\left(|A_{i,j}|+\sum_{Q\in F_{i,j}}|Q|\right)\\
&\leq&(C+1)\sum_{(i,j)}|A_{i,j}|=(C+1)|S|
\end{eqnarray*}
where $C$ is the consant in (\ref{ddjpe.eq:fijaij}).

\medskip

\noindent{\sl Proof of Property \ref{ddjpe.pointii}.}

\smallskip

If $(\eta,a)\in S$, that is if $(\eta,a)\in A_{i,j}$ for some $(i,j)$ then there exists $(\eta',a')\in A_{i,j}$
such that $U_\nu(\eta',a')>\lambda$. But then 
$|\eta^{-1}\eta'|\leq 2^{i_0-i-2}\leq \frac{a'}{2}$ so that,
for $\xi\in\fN$,
\begin{eqnarray*}
a'+|\xi^{-1}\eta'|&\geq&a'-|\eta^{-1}\eta'|+|\xi^{-1}\eta|\geq
\frac{a'}{2}+|\xi^{-1}\eta|\\
&\geq&\frac{1}{4}(a+|\xi^{-1}\eta|).
\end{eqnarray*}
From this, we immediately get that
$$
\alpha<U_\nu(\eta',a')\leq CU_\nu(\eta,a).
$$

\medskip

\noindent{\sl Proof of Property \ref{ddjpe.pointiii}.}

\smallskip

Set $S_i=\dst\bigcup_{j}A_{i,j}$ the set of authorized pieces of order $i$ and write
$$
S_i=T_i\times[2^{i_0-i-1},2^{i_0-i}].
$$
Set
$$
U_i(\eta)=\int_{S_i}\frac{a^{\Gamma-1}}{(a+|\xi^{-1}\eta|)^{Q+\Gamma}}
\,\mbox{d}\lambda(\xi)\,\mbox{d}a
$$
for the part of $U_S$ issued from pieces of order $i$.

\begin{lemma}
\label{ddjpe.lem:sjogrenlem2}
There exists a constant $C_2$ such that, for all $\eta\in\fN$ and all $i\geq0$, $U_i(\eta)\leq C_2$. Moreover,
if $p\geq -i$ and $\mathrm{d}(\eta,T_i)>2^{i_0+p}$, then $U_i(\eta)\leq C_22^{-(p+i)\Gamma}$.
\end{lemma}

\begin{proof}[Proof of Lemma \ref{ddjpe.lem:sjogrenlem2}]
By definition
\begin{equation}
U_i(\eta)=\int_{2^{i_0-i-1}}^{2^{i_0-i}}a^{\Gamma-1}\int_{T_i}
\frac{\mbox{d}\lambda(\xi)}{(a+|\xi^{-1}\eta|)^{Q+\Gamma}}\,\mbox{d}a
\label{ddjpe.eq:sjogrenl2}
\end{equation}
Thus
\begin{eqnarray*}
U_i(\eta)&=&\int_{2^{i_0-i-1}}^{2^{i_0-i}}\frac{1}{a^{Q+1}}
\int_{T_i}\frac{\mbox{d}\lambda(\xi)}{(1+|\xi^{-1}\eta|/a)^{Q+\Gamma}}
\,\mbox{d}a\\
&=&\int_{2^{i_0-i-1}}^{2^{i_0-i}}\frac{1}{a}\int_{T_i/a}
\frac{\mbox{d}\lambda(\zeta)}{(1+|\zeta^{-1}\cdot(\eta/a)|)^{Q+\Gamma}}
\,\mbox{d}a,
\end{eqnarray*}
changing variable $\zeta=\xi/a$. By translation invariance of
$\mbox{d}\lambda$ we thus get that
$$
U_i(\eta)\leq\int_{2^{i_0-i-1}}^{2^{i_0-i}}\frac{1}{a}\int_{\fN}
\frac{\mbox{d}\lambda(\zeta)}{(1+|\zeta|)^{Q+\Gamma}}\,\mbox{d}a
\leq C_2.
$$

Further, if $\mathrm{d}(\eta,T_i)>2^{i_0+p}$ then, for $\xi\in T_i$,
$|\xi^{-1}\eta|\geq2^{i_0+p}$ so that,
from (\ref{ddjpe.eq:sjogrenl2}), we deduce that
\begin{eqnarray*}
U_i(\eta)&\leq&2^{(i_0-i)\Gamma}\int_{T_i}
\frac{\mbox{d}\lambda(\xi)}{(2^{i_0-i-1}+|\xi^{-1}\eta|)^{Q+\Gamma}}\\
&\leq&
2^{(i_0-i)\Gamma}\int_{|\xi^{-1}\eta|\geq 2^{i_0+p}}
\frac{\mbox{d}\lambda(\xi)}{|\xi^{-1}\eta|^{Q+\Gamma}}\\
&\leq&2^{(i_0-i)\Gamma}\int_{|\zeta|>2^{i_0+p}}
\frac{\mbox{d}\lambda(\zeta)}{|\zeta|^{Q+\Gamma}}\\
&=&C2^{-(i+p)\Gamma}
\end{eqnarray*}
when integrating in polar coordinates.
\end{proof}

Now, for every $\eta\in\fN$, there exists an $m$
such that
\begin{equation}
\label{ddjpe.eq:sjo*}
C_22^{-m-1}<U_j(\eta)\leq C_22^{-m}
\end{equation}
(where $C_2$ is the constant of Lemma \ref{ddjpe.lem:sjogrenlem2}). From Lemma \ref{ddjpe.lem:sjogrenlem2},
we get that $\mathrm{d}(\eta,T_j)<2^{i_0+(m+1)/\Gamma-j}$.

On the other hand, by construction, if $i\leq j$ then
$$
\mathrm{d}(T_i,T_j)>2^{i_0-i+\frac{j-i}{Q+\Gamma}+1}
$$
so that
$$
\mathrm{d}(\eta,T_i)>\mathrm{d}(T_i,T_j)-\mathrm{d}(\eta,T_j)
>2^{i_0-i+\frac{j-i}{Q+\Gamma}+1}-2^{i_0+(m+1)/\Gamma-j}.
$$
It follows that $\mathrm{d}(\eta,T_i)>
2.2^{i_0-i+\frac{j-i}{Q+\Gamma}}-2^{i_0+(m+1)/\Gamma-j}$.
Further, if $i_0+(m+1)/\Gamma-j<i_0-i+\frac{j-i}{Q+\Gamma}$,
in particular, if $j-i\geq (m+1)/\Gamma$, then
$$
\mathrm{d}(\eta,T_i)> 2^{i_0-i+\frac{j-i}{Q+\Gamma}}.
$$
From Lemma \ref{ddjpe.lem:sjogrenlem2} we then get that
\begin{equation}
\label{ddjpe.eq:abovecomp}
U_i(\eta)\leq C_22^{-\frac{\Gamma}{Q+\Gamma}(j-i)}
\end{equation}
for $i\leq j-m$.

\medskip

We will now prove by induction on $j$ that
there exists a constant $C_2$ for which,
{\sl for every $\eta>0$ and every
$j\geq0$, there exists a permutation $\sigma=\sigma_{j,\eta}$
of $\{0,\ldots,j\}$ such that,
for each $i\in\{0,\ldots,j\}$, $U_i(\eta)\leq C_22^{-\frac{\Gamma}{Q+\Gamma}\sigma(i)}$.}
\smallskip

It then immediatly follows that $\sum U_i$ is convergent and uniformly bounded as desired.

\medskip

For $j=0$, this is just Lemma \ref{ddjpe.lem:sjogrenlem2}. Assume now the hypothesis is true up to order $j-1$.

Let $\eta\in\fN$ and let $m$ be such that $C_22^{-m-1}\leq U_j(\eta)\leq C_22^{-m}$.

\smallskip

\noindent --- If $m\geq j$, then $U_j(\eta)\leq C_22^{-\frac{\Gamma}{Q+\Gamma}j}$.
Further, by induction hypothesis, there exists a permutation $\sigma_{j-1,\eta}$ of $\{0,\ldots,j-1\}$
such that, for $i\in\{0,\ldots,j-1\}$, $U_i\leq 2^{-\frac{\Gamma}{Q+\Gamma}\sigma_{j-1,\eta}(i)}$. 
It is then enough to extend $\sigma_{j-1,\eta}$ by setting
$\sigma_{j,\eta}(i)=\sigma_{j-1,\eta}(i)$ if $i<j$ and
$\sigma_{j,\eta}(j)=j$.

\smallskip

\noindent --- Otherwise, $m< j$ and the \eqref{ddjpe.eq:abovecomp}
shows that, for $i=0,\ldots,j-m-1$,
$U_i(\eta)\leq C_22^{-\frac{\Gamma}{Q+\Gamma}(j-i)}$.

By the induction hypothesis, $U_{j-m}(\eta),\ldots,U_{j-1}(\eta)$
are bounded by $m-1$ different elements of
$\bigl\{C_2,C_22^{-\frac{\Gamma}{Q+\Gamma}},\ldots,
C_22^{-\frac{\Gamma}{Q+\Gamma}(j-1)}\bigr\}$.
But these are decreasing, so we may as well assume
that they are bounded by the $m-1$ first elements of the family.
In other words, there exists a one-to-one mapping $\sigma_1$ from
$\{j-m,\ldots,j-1\}$ to $\{0,\ldots,m-1\}$ such that,
for $i=j-m,\ldots,j-1$, $U_{i}(\eta)\leq C_22^{-\frac{\Gamma}{Q+\Gamma}\sigma(i)}$.

Finally, as $U_j(\eta)\leq C_22^{-\frac{\Gamma}{Q+\Gamma}m}$,
if we set
$$
\sigma_{j,\eta}(i)=\begin{cases} j-i&\mbox{if }i=0,\ldots,j-m-1\\
\sigma_1(i)&\mbox{if }i=j-m,\ldots,j-1\\
m&\mbox{if }i=j
\end{cases}.
$$
the proof of the induction is completed.
\end{proof}
We have thus completed the proof \eqref{ddjpe.eq:sjogren2}.
\end{proof}

\end{document}